%% file: main.tex
\documentclass{amsart}
\usepackage{import}
\usepackage{./preamble/Pakete}
\usepackage{./preamble/Befehle}
\bibliography{bibfinal.bib}

\begin{document}

\import{Inhalt/}{title}

\import{Inhalt/}{Intro.tex}
\import{Inhalt/}{setup.tex}
\import{Inhalt/}{singularKtheoryBimodules.tex}
\import{Inhalt/}{ktheoreticsoergel2.tex}

\printbibliography

\end{document}

%% file: Inhalt/title.tex
\title{The Parabolic K-motivic Hecke Category}

\author{Jens Niklas Eberhardt}
\email{mail@jenseberhardt.com}
\address{Johannes Gutenberg-Universität Mainz, Institut für Mathematik, Staudingerweg 9, 55128 Mainz, Germany}

\author{Arnaud Eteve}
\email{eteve@mpim-bonn.mpg.de}
\address{Max Planck Institute for Mathematics, Vivatsgasse 7, 53111 Bonn, Germany}

\begin{abstract}
We define and study the parabolic $K$-motivic Hecke category of a (possibly disconnected) Kac--Moody group. Our main result is a combinatorial description via singular $K$-theory Soergel bimodules which arise from the equivariant algebraic $K$-theory of parabolic Bott--Samelson resolutions. 
In the spherical affine case, the $K$-motivic Hecke category serves as one side of a conjectural quantum $K$-theoretic derived Satake equivalence, addressing a conjecture of Cautis--Kamnitzer.
\end{abstract}
\maketitle

%% file: Inhalt/Intro.tex
\section{Introduction}
Let $G\supset B\supset T$ be a (possibly disconnected) Kac--Moody group with a fixed Borel subgroup and maximal torus. In this article, we will define the \emph{parabolic $K$-motivic Hecke category}
$$\Hecke_{P,Q}\subset \DK(P\bs G/Q)$$
attached to a pair of standard parabolic subgroups of finite type $B\subset P,Q\subset G$. Our main result is a `Soergel-type' description of this category that depends only on the action of parabolic Weyl groups on the representation ring $\RR=\Z[X(T)]$.
\begin{theorem*}[\Cref{thm:K-is-an-equivalence}]
    Assume that the simple roots of $G$ are linearly independent and $\pi_1(G)$ is free. Then there is an equivalence
    \[\begin{tikzcd}
        {\Hecke_{P,Q}} & {\Ch^b(\SBim_{P,Q})}
        \arrow["\sim", from=1-1, to=1-2]
    \end{tikzcd}\]
    with the category of chain complexes of singular $K$-theory Soergel bimodules.
\end{theorem*}
We refer to \cite{eberhardtKMotivesKoszulDuality2022, eberhardtKtheorySoergelBimodules2024,eberhardtUniversalKoszulDuality2025a} which treat the case of $P=Q=B$ and connected groups $G$ as well as the work of Soergel--Virk--Wendt in a similar motivic setting \cite{soergelEquivariantMotivesGeometric2018}.

In the rest of the introduction, we will explain the various ingredients of this equivalence as well as the relation to a quantum $K$-theoretic Satake equivalence that was one main motivation for this work.
\subsection{Genuine equivariant $K$-motives}
The category $\Hecke_{P,Q}$ arises within the formalism of \emph{reduced genuine equivariant $K$-motives} $\DK$, which is based on Hoyois' genuine equivariant stable motivic homotopy category \cite{hoyoisSixOperationsEquivariant2017,hoyoisCdhDescentEquivariant2020}. This formalism provides a categorical version of genuine equivariant algebraic $K$-theory and is equipped with a six-functor formalism. We refer to \Cref{sec:kmotives} for details and only highlight the most important features here.

The category $\DK(\Xx)$ is defined for `nice enough' stacks $\Xx$ and shares many formal properties with the derived categories of $D$-modules and constructible sheaves. In particular there are the six functors $f^*,f_*,f_!,f^!$ (for $f$ representable) and $\otimes, \iHom$ satisfying the usual compatibilities. 

A key distinguishing feature is that $\DK$ computes \emph{genuine equivariant algebraic $K$-theory}. For example, for a reductive group $G$ with maximal torus $T$ and Weyl group $W$, there is an equivalence
$$\DKT(\pt/G)\simeq \op{Perf}(\RR_G)$$
between Tate $K$-motives (those stably generated by the unit object) and perfect complexes over the representation ring $\RR_G=K_0(\pt/G)$.\footnote{We note that the topological version of $K$-motives, namely modules over the complex $K$-theory spectrum $KU$, is $[2]$-periodic and yields $[2]$-periodic perfect complexes here. This category does not admit a non-trivial $t$-structure or weight structure which we want to avoid.}
\subsection{$K$-motivic Hecke category}
We now introduce the main object of study.
\begin{definition*}[see \Cref{def:heckecat}]
The \emph{parabolic $K$-motivic Hecke category} is the full stable thick subcategory
$$\Hecke_{P,Q}=\langle \Delta_w \rangle_{\op{stable},\op{thick}}\subset \DK(P\bs G/Q)$$
generated by the standard objects $\Delta_w=i_{w,!}\un$ corresponding to the different Bruhat strata of the parabolic Hecke stack $P\bs G/Q$.
\end{definition*}

The categories $\Hecke_{P,Q}$ for all pairs $P,Q\subset G$ carry a natural convolution structure: there are bifunctors $\Hecke_{P,Q}\times \Hecke_{Q,R}\to \Hecke_{P,R}$. In the case $P=Q$, this endows $\Hecke_{P,P}$ with a monoidal structure, and the category admits a self-duality exchanging standard and costandard objects.

We will often be able to reduce statements to the Borel-case $\Hecke_{B,B}$ for a connected group $G$ which we studied in previous work \cite{eberhardtUniversalKoszulDuality2025a}. 
For this, we use the map  $\pi:B\bs G/B\to P\bs G/Q$ and show that it satisfies an equivariant version of the projective bundle formula (\Cref{cor:projectivebundleformula}) using the Pittie--Steinberg theorem \cite{steinbergTheoremPittie1975}.

Within $\Hecke_{P,Q}$ we consider \emph{pure} objects $\Hecke_{P,Q}^{\pure}$ generated under direct sums and summands by pushforwards of Bott--Samelson resolutions. The reader might compare this definition with Soergel's category $\mathcal{K}$ or parity sheaves, see \cite{soergelRelationIntersectionCohomology2000a,juteauParitySheaves2014}.

We show that pure objects satisfy a \emph{pointwise purity} condition (\Cref{thm:pointwisepuritybottsamelson}) which implies that their mapping spaces are discrete. This implies the following formality theorem.
\begin{theorem*}[\Cref{corol:weight-complex-functor-is-an-equiv}]
    There is an equivalence of categories
    $$\Hecke_{P,Q}\simeq \Ch^b(\Hecke_{P,Q}^{\pure})$$
    between the parabolic $K$-motivic Hecke category and the category of chain complexes of pure objects.
\end{theorem*}
\subsection{Singular $K$-theory Soergel bimodules}
We now describe the combinatorial side of the main equivalence, which is a $K$-theoretic version of singular Soergel bimodules \cite{williamsonSingularSoergelBimodules2010}. The category $\SBim_{P,Q}$ of \emph{singular $K$-theory Soergel bimodules} is defined as the idempotent completion of the category of Bott--Samelson bimodules. For a connected component $G^\circ$, a Bott--Samelson bimodule is a finite iterated tensor product of the form
$$\RR_P\otimes_{\RR_{Q_1}}\RR_{P_1}\otimes_{\RR_{Q_2}}\cdots \otimes_{\RR_{Q_n}}\RR_Q$$
where $P,P_i,Q_i,Q$ are parabolic subgroups. This construction generalizes to disconnected groups by decorating the tensor products with elements of $\Gamma=\pi_0(G)$ that twist the bimodule structures.

To relate $\SBim_{P,Q}$ to $\Hecke_{P,Q}^{\pure}$, we consider the functor $$\Kyp:\Hecke_{P,Q}\to \Dd(\RR_P\otimes\RR_Q)$$ that sends a $K$-motive to its $K$-theory as a bimodule. We prove that $\Kyp$ is monoidal (see \Cref{thm:Strict-Monoidality}) and establish an analog of Soergel's Erweiterungssatz. 
\begin{theorem*}[see \Cref{thm:K-is-an-equivalence}]
   The functor $\Kyp$ induces an equivalence 
    \[\begin{tikzcd}
        {\Hecke_{P,Q}^{\pure}} & {\SBim_{P,Q}.}
        \arrow["\sim", from=1-1, to=1-2]
    \end{tikzcd}\]
\end{theorem*}
\subsection{The spherical case} 
A main motivation for the current work is that there should be an equivalence between the \emph{spherical affine} $K$-motivic Hecke category and a certain category of representations of a quantum group.
\begin{remark}
    Such an equivalence is analogous to the equivalence between equivariant Springer $K$-motives on the nilpotent cone and representations of the affine Hecke algebra, see \cite{kazhdanProofDeligneLanglandsConjecture1987,eberhardtKmotivesSpringerTheory2024}.
\end{remark}
The current work yields a combinatorial description of the $K$-motivic side and should hence be seen as the first half of a bridge under construction. We will now explain some more details.

Let $\mathring{G}$ be a complex reductive group with a maximal torus $\mathring{T}$ and Weyl group $W_f$. We assume that $\mathring{G}$ has a simply connected derived subgroup. Then we may consider the semidirect product of the loop group and the group of loop rotations
$$G=L\mathring{G}\rtimes \Gm$$ which is a (possibly disconnected) Kac--Moody group where
$\pi_0(G)=\pi_1(\mathring{G})$.
The positive loop group yields a maximal parabolic subgroup of finite type
$$P=L^+\mathring G\rtimes \Gm\subset G$$
and we may hence consider the double quotient
$$ P\bs G/P=(L^+\mathring G\rtimes \Gm)\backslash \op{Gr}_{\mathring G}$$
which is a quotient of the affine Grassmannian of $\mathring G$.
We obtain the \emph{$K$-motivic spherical affine Hecke category}
$$\Hecke_{P,P}\subset \DK(P\bs G/P)=\DK((L^+\mathring G\rtimes \Gm)\backslash \op{Gr}_{\mathring G})$$
which is naturally a bimodule over the representation ring 
$$\RR_P=\RR_{\mathring{G}}[q^{\pm 1}]=\Z[X(\mathring{T})]^{W_f}[q^{\pm 1}].$$
Now \Cref{corol:weight-complex-functor-is-an-equiv} provides a combinatorial description of this category
in terms of chain complexes of singular $K$-theory Soergel bimodules over $\RR_P.$ As a first hint that this is related to the quantum group, we observe that the quantum Harish-Chandra isomorphism \cite{rossoAnaloguesFormeKilling1990} yields an isomorphism 
$$Z(U_q(\mathring{\mathfrak{g}}))\cong \RR_P\otimes_{\Z[q^{\pm 1}]}\mathbb{C}(q).$$
In the next section, we discuss a much stronger hint.
\subsection{Quantum Satake for $\op{SL}_n$}
In the following, we scalar extend all categories along $\Z[q^{\pm 1}]\subset \mathbb C(q)$.

Cautis--Kamnitzer construct a so-called quantum $K$-theoretic Satake equivalence for $\mathring G=\op{SL}_n$, see \cite{cautisQuantumKtheoreticGeometric2018}. The relation with the current work can be summarized in the following diagram
\[\begin{tikzcd}
	{\Hh_{P,P}^{\pure}} & {\SBim_{P,P}} \\
	{\op{KConv}^{\mathring{G}\times \Gm}(\op{Gr}_{\mathring G})} & {\mathcal{AS}p_n(q)}.
	\arrow["\sim", from=1-1, to=1-2]
	\arrow["\wr", from=2-1, to=1-1]
	\arrow["\wr", dashed, from=2-2, to=1-2]
	\arrow["\sim", from=2-1, to=2-2]
\end{tikzcd}\]
Here, the bottom left category is the so-called \emph{$K$-theoretic convolution category} considered by Cautis--Kamnitzer. Hom-spaces in this category are given by $\op{G}$-theory groups $\op{G}_0(\Xx_1\times_{P\bs G/P}\Xx_2)$ where $\Xx_i$ denote certain resolutions of Schubert varieties, and composition is given via convolution. This category is equivalent to the category of pure objects in the spherical affine $K$-motivic Hecke category. We refer to a similar comparison between pure $K$-motives on $B\bs G/B$ and a convolution category in \cite[Remark 4.5 (3)]{eberhardtKtheorySoergelBimodules2024} and also the proof of \Cref{thm:pushforwardpittiesteinberg}.

The bottom arrow is Cautis--Kamnitzer's diagrammatic description of the $K$-theoretic convolution category in terms of the category of annular spiders $\mathcal{AS}p_n(q)$. Combining this with the results of our paper, we see that there is an equivalence between singular $K$-theory Soergel bimodules and the category $\mathcal{AS}p_n(q)$ which might be of independent interest.
See also \cite{eliasQuantumSatakeType2017} for similar equivalences for classical singular Soergel bimodules. 

\begin{remark}
    While the top arrow in the above diagram works for any type, the bottom arrow in the diagram relies on type $A$. Objects in the Cautis--Kamnitzer $K$-theoretic convolution category correspond to products of miniscule Schubert varieties in the affine Grassmannian, which only provide resolutions of all Schubert varieties if one is in type $A$. However, one should slightly modify the definition of the category, such that the objects correspond to general Bott--Samelson varieties. In this case, the left arrow is still an equivalence. Composing the left and top arrows then yields a combinatorial description of the $K$-theoretic convolution algebra, answering a question of Cautis--Kamnitzer.
\end{remark}

In the same work, Cautis--Kamnitzer also construct an equivalence 
\[\begin{tikzcd}
	{\mathcal{AS}p_n(q)} & {\mathcal{O}_q(\mathring{G}/\op{Ad}(\mathring{G}))^{\op{min}}\text{-mod}}
	\arrow["\sim", from=1-1, to=1-2]
\end{tikzcd}\]
where the right hand side denotes the additive category of $U_q(\mathring{\mathfrak{g}})$-equivariant $\mathcal{O}_q(\mathring{G})$-modules generated by modules of the form $\mathcal{O}_q(\mathring{G})\otimes L(\lambda_1)\otimes\dots L(\lambda_m)$.
Passing to homotopy categories, we obtain an equivalence
\[\begin{tikzcd}
	\Hh_{P,P} & \Ch^b(\mathcal{O}_q(\mathring{G}/\op{Ad}(\mathring{G}))^{\op{min}}\text{-mod})=:\op{Perf}_q^{\op{min}}(\mathring{G}/\op{Ad}(\mathring{G}))
	\arrow["\sim", from=1-1, to=1-2]
\end{tikzcd}\]
with a $q$-deformed category of perfect complexes on the adjoint quotient. It would be very desirable to extend this to all types. 

\subsection{Further directions}
There should be an equivalence between the parabolic $K$-motivic Hecke category and bi-Whittaker monodromic sheaves for the Langlands dual group generalizing the universal Koszul duality considered in \cite{eberhardtKtheorySoergelBimodules2024,taylorUniversalMonodromicTilting2023,eberhardtUniversalKoszulDuality2025a}. We also refer to the recent work \cite{chenLanglandsDualRealization2023a, dhillonEndoscopyMetaplecticAffine2025, dhillonTameLocalBetti2025}.

In a different direction, it would be interesting to compare the $K$-motivic Hecke category to categories of coherent sheaves via the categorical Chern character \cite{hoyoisHigherTracesNoncommutative2017} which yields a close relation to recently proposed potent categorical representations \cite{ben-zviPotentCategoricalRepresentations2025}.

\subsection{Acknowledgments}
The first author is grateful to Joel Kamnitzer for pointing us to the quantum $K$-theoretic Satake. We warmly thank Alexander Bravermann, Ilya Dumanski and Antoine Labelle for interesting discussions on the derived Satake equivalence. The first author was supported by Deutsche Forschungsgemeinschaft (DFG), project number 45744154, Equivariant K-motives and Koszul duality. The second author was supported by the Max Planck Institut for Mathematics in Bonn during the preparation of this paper.

\subsection{Notations and conventions}

\subsubsection{Group and root theoretic data}

We let $G,B,T,U$ be Kac-Moody groups over $k$, an algebraically closed field of characteristic $0$, with a chosen Borel subgroup, and we let  $(X^*, \Phi, \Phi^+, \Delta, X_*, \Phi^{\vee}, \Phi^{\vee, +}, \Delta^{\vee})$ be the usual sets of (co)characters of $T$, (co)roots, positive (co)roots and simple (co)roots. In most sections we assume that $G_{\mathrm{der}}$ is simply connected; this assumption will be recalled as needed. More details are given in Section \ref{sec:Kac-Moody-groups}. In Sections \ref{sec:Singular-K-theory-SBim} and \ref{sec:Main-section}, the group $G$ is allowed to be disconnected; we refer to Section \ref{sec:disconnected-groups} for the precise assumptions on $G$. 

\subsubsection{Categories}

We shall use the language of $\infty$-categories of Lurie \cite{lurieHigherToposTheory2009}. For an $\infty$-category, $\Map$ denotes the mapping space and $\Hom = \pi_0(\Map)$. For categories, we denote by $\mathrm{Pr}$ the category of presentable stable $\infty$-categories with continuous functors and by $\mathrm{Pr}^{\omega}$ the (not full) subcategory of presentable stable compactly generated $\infty$-categories with continuous compact preserving functors. 

For a (classical) $1$-category, $\Hom$ denotes the usual $\Hom$ space. 

\subsubsection{Algebras and modules}

Let $R$ be a ring. We shall denote by $R\text{-Mod}$ the abelian category of $R$-modules, by $\Dd(R)$ the full $\infty$-derived category of $R\text{-Mod}$, and by $\Perf(R) \subset \Dd(R)$ the full subcategory of perfect complexes (i.e., complexes which are isomorphic to finite complexes of finitely generated projective $R$-modules).

\subsubsection{Morita category}\label{sec:Morita-category}

For a symmetric monoidal category $\mathcal{C}$ (either $1$ or $\infty$), we denote by $\Alg(\mathcal{C})$ the category of associative algebras in $\mathcal{C}$. Given an algebra $A$ and a right (resp. left) $A$-module $M$ (resp. $N$), if it is defined, we denote by $M \otimes_A N$ the relative tensor product. The main example of such a $1$-category is $k\text{-Mod}$ for some commutative ring $k$. At the $\infty$-categorical level, the main example is the category $\Pr$ or $\Pr_{k}$ of $k$-linear presentable categories equipped with the Lurie tensor product. Let $\mathcal{C}$ be a category where the relative tensor product exists (note that by \cite{higherAlgebra} this always exists provided geometric realization exists); then we will write $\Morita(\mathcal{C})$ for the category whose objects are algebras and whose morphisms are bimodules, and composition is given by relative tensor of bimodules. In the $1$-categorical case, this category is defined unambiguously; for the $\infty$-categorical setting, we refer to \cite[Section 4.4.3]{higherAlgebra}. We remark that $\Morita(\mathcal{C})$ should really be a $2$ (or $(\infty, 2)$)-category; however, we will view it only as a $1$ (or $(\infty, 1)$)-category by removing the non-invertible $2$-morphisms. 
\begin{remark}\label{rmk:Grothendieck-construction}
In this paper, we will consider functors $I \to \Morita(\mathcal{C})$ from some index category $I$. It is common in the literature to use the language of monoidal $2$-categories; we will not use this language. This is mainly to avoid having to define $(\infty, 2)$-categories, but it should be noted that the Grothendieck construction gives a dictionary between monoidal $2$-categories $\mathcal{D}$ fibered over $I$ and functors $I \to \Morita(\mathcal{C})$ for some appropriate symmetric monoidal category $\mathcal{C}$. 
\end{remark}

%% file: Inhalt/setup.tex
\section{Categories of $K$-motives} \label{sec:kmotives}

We recall important definitions and properties of $K$-motives, following \cite{eberhardtUniversalKoszulDuality2025a,eberhardtKmotivesSpringerTheory2024}.
Throughout the document, let $\point=\Spec(k)$
for an algebraically closed field $k$ of characteristic zero. Recall that a stack $\Xx/k$ is \emph{\nice} if it can be represented as a quotient stack $X/G$ where $G/k$ is linearly reductive, and $X/k$ is of finite type and $G$-quasi-projective, so admits an equivariant embedding into $\mathbb P(V)$ for a representation $V$.
\subsection{Reduced $K$-motives}\label{sec:reducedkmotives}
The work of Hoyois \cite{hoyoisSixOperationsEquivariant2017, hoyoisCdhDescentEquivariant2020} defines a \emph{genuine equivariant stable motivic homotopy category} $\SH(\Xx)$ over $\Xx$ and shows that there is a Cartesian ring spectrum $\KGL_\Xx$ in $\SH(\Xx)$ representing $\A^1$-homotopy invariant algebraic $K$-theory. Using this, one may define a category of \emph{genuine equivariant $K$-motives} on $\Xx$ (or just $K$-motives for short) as the module category
$$\DKbig(\Xx)=\Mod_{\KGL_\Xx}(\SH(\Xx)).$$
In \emph{loc.~cit.~} it is further shown that there is a functor between the category of correspondences of linearly reductive stacks with quasi-projective morphisms
\begin{align*}
    \op{Corr}(\op{Stk})&\to \op{Pr}: \Xx \mapsto \DKbig(\Xx),\, (\Xx\xleftarrow{f} \Zz \xrightarrow{g} \Yy)\mapsto g_!f^*
\end{align*}
which yields a six functor formalism with $f^*,f_*,f_!,f^!$ for quasi-projective morphisms and bifunctors $\iHom,\otimes.$ We will denote the tensor unit (or constant $K$-motive) by $\un=\un_\Xx\in \DKbig(\Xx).$

\begin{remark}\label{remark:comparison-pull-back}
A remarkable property of $K$-motives is \emph{Bott-periodicity} which implies that $(1)[2]$ is isomorphic to the identity functor. This has the following interesting consequences.
Let $f : \Xx \to \Yy$ be a smooth quasi-projective morphism. Then there is a natural isomorphism $f^* = f^!$. If $f$ is an affine bundle, then there are isomorphisms $f_!\un\cong f_*\un=\un$.
\end{remark} 
We will work with a slight modification of the category of $K$-motives in which we remove the higher algebraic $K$-groups $K_{>0}(\pt)$ of the ground field. For this, we define the category of 
\emph{reduced $K$-motives}, see also \cite{eberhardtIntegralMotivicSheaves2023}, as the tensor product
$$\DK(\Xx)=\DKbig(\Xx)\otimes_{\op{K}(\point)}\op{K}_{0}(\point).$$ Here $\op{K}(\point)$ acts on $\DKbig(\Xx)$ via the maps $$\op{K}(\point)\to \op{K}(\pt/G)=\Map_{\DKbig(\pt/G)}(\un,\un)\stackrel{f^*}{\to} \Map_{\DKbig(\Xx)}(\un,\un)$$
where $\Xx=X/G$ is a presentation of the linearly reductive stack $\Xx$. The natural map $\DKbig\to \DK$ is then compatible with the six functors, and mapping spaces in $\DK$ are related to \emph{reduced $K$-theory} and \emph{reduced $G$-theory}
\begin{align*}
    \Map_{\DK(\Xx)}(\un,\un) &\simeq \op{K}(\Xx)_\red = \op{K}(\Xx)\otimes_{\op{K}(\point)}\op{K}_{0}(\point) \text{ and}\\
    \Map_{\DK(\Yy)}(\un,f^!\un)& \simeq \op{G}(\Yy)_\red = \op{G}(\Yy)\otimes_{\op{K}(\point)}\op{K}_{0}(\point),
\end{align*}
for $\Xx/k$ smooth and $f:\Yy\to\Xx$ a quasi-projective map of \nice\ stacks. The computation further simplifies for stacks that admit a stratification into affine bundles over classifying stacks of linear algebraic groups, in which case $\op{G}(\Yy)_\red=G_0(\Yy)$, see \cite{eberhardtKmotivesSpringerTheory2024}.

Lastly, we will need an extension of the formalism to infinite-dimensional stacks. Indeed, the formalism $\DK$ can be extended to linearly reductive ind-pro stacks $\Xx$. These are stacks that admit a presentation
$\Xx=\colim_n \lim_m \Xx_{n,m}$ by linearly reductive $\Xx_{n,m}$, such that maps in the $n$ and $m$-direction are affine bundles and closed immersions, respectively. For more details, we refer to \cite[Section 3.4]{eberhardtUniversalKoszulDuality2025a}.

\subsection{$K$-motives on classifying stacks}
We now discuss the special case of reduced $K$-motives on classifying stacks. For this, we define for a stack $\Xx$ the category of Tate reduced $K$-motives
\begin{align*}
	\DKT(\Xx)=\langle \un \rangle_{\op{stable},\op{thick}}\subset \DK(\Xx)
\end{align*}
as the full stable thick subcategory generated by the unit object. 
\begin{remark}
	Usually, one would take all Tate objects $\un(n)$ for $n\in\Z$ as generators. However, for $K$-motives $\un(n)=\un[-2n]$ by Bott periodicity, see \Cref{remark:comparison-pull-back}.
\end{remark}
An important statement that we simply refer to as `homotopy invariance' will often allow us to pass to the case of reductive groups.
\begin{lemma}\label{lem:homotopyinvariance}
	Let $P=L\ltimes U$ be a semidirect product of a linear algebraic group with a pro-unipotent group with a presentation $U=\lim U/U(n)$, where the subgroups $U(n)$ are normal in $P$. Moreover, assume that for each $n$ the quotient $U(n)/U(n+1)$ is $P$-equivariantly isomorphic to a representation of $P$. 
	Denote by $\pi: \pt/L\to \pt/P$ the projection.
	Then the functor $\pi^*: \DK(\pt/P)\to \DK(\pt/L)$ is fully faithful.
\end{lemma}
\begin{proof}
	This is \cite[Lemma 3.1]{eberhardtUniversalKoszulDuality2025a}
\end{proof}
\begin{remark}
	The fully-faithfulness implies that $\pi_*\pi^*=\id$ and that the functor $\pi^*$ induces an equivalence on the category of Tate motives.
\end{remark}
\begin{remark}
	The assumptions of \Cref{lem:homotopyinvariance} are, for example, fulfilled for $P=L\ltimes U$ a finite-type parabolic subgroup in a Kac--Moody group, where $L$ is the Levi-factor and $U$ the unipotent radical. The subgroups $U(n)$ are generated from all root subgroups whose $P$-relative height is $\geq n.$
\end{remark}
The computation of mapping spectra for reduced $K$-motives in \Cref{sec:reducedkmotives} yields the following description for the Tate subcategory.
\begin{lemma}\label{lem:DescriptionOfDKT}
Let $G$ be a reductive group. Then there is an equivalence
 $$\RR_G=K_0(\pt/G)\simeq \Map_{\DK(\pt/G)}(\un,\un)$$
between the (discrete) representation ring of $G$ and the mapping endomorphism spectrum of $\un\in \DK(\pt/G)$. Moreover, the functor $\Map(\un,-)$ yields an equivalence
\begin{align*}
	\DKT(\pt/G) \stackrel{\sim}{\to} \op{Perf}(\RR_G)
\end{align*}
between Tate $K$-motives on $\pt/G$ and the category of perfect complexes over $\RR_G$.
\end{lemma}
The representation ring $\RR_G$ admits a concrete description. Namely, denote by $T\subset G$ a maximal torus and by $W=N_G(T)/T$ the Weyl group (which is a discrete finite group). Moreover, denote by $X(T)$ the character lattice of $T$. Then the Chevalley restriction theorem shows that the restriction map yields an isomorphism
$$\RR_G=\RR_T^W=\Z[X(T)]^W$$
with the Weyl group invariants of the representation ring of $T$.

Our next goal is to study the functoriality of \Cref{lem:DescriptionOfDKT}.
\begin{theorem}\label{thm:pushforwardpittiesteinberg}
Let $G$ be a reductive group with a simply-connected derived subgroup and $B$ be a Borel subgroup and $W$ the Weyl group. Denote by $p: \pt/B\to\pt/G$ the projection. Then $p_!(\un)=p_*(\un)=\un^{\oplus |W|}\in \DK(\pt/G)$.
\end{theorem}
\begin{proof}
	Since the map $\DKbig\to \DK$ is compatible with the six functors, we might as well prove the result in $\DKbig(\pt/G)$. Hence, in the following, all hom-spaces are in $\DKbig(\pt/G)$. 
	
	To follow the exposition in \cite{kazhdanProofDeligneLanglandsConjecture1987} more closely, we rewrite our map $p$ in the form $p:G\bs G/B\to G\bs \pt.$
	For proper maps of linearly reductive stacks $p_i:\Xx_i\to \Yy$ such that $\Xx_i$ is smooth, we have the general formula
	$$\Hom_{\DKbig(\Yy)}(p_{1,*}\un,p_{2,*}\un)=G_0(\Xx_1\times_{\Yy}\Xx_2)$$
	such that composition of homomorphisms corresponds to convolution on $G$-theory, see \cite{eberhardtKmotivesSpringerTheory2024}. We hence obtain
	\begin{align*}
		\End(p_*\un)&=G_0(G\bs(G/B\times G/B)), \\
		\Hom(\un,p_*\un)&=G_0(G\bs G/B)=\RR_B, \\
		\Hom(p_*\un,\un)&=G_0(G\bs G/B)=\RR_B \text{ and }\\
		\End(\un)&=G_0(G\bs \pt)=\RR_G.
	\end{align*}
	In one direction, composition of homomorphisms can be described via the diagram
	
	\[\begin{tikzcd}
		{\Hom(\un,p_*\un)\times \Hom(p_*\un,\un)} & {\End(\un)} \\
		{G_0(G\bs G/B)\times G_0(G\bs G/B)} & {G_0(G\bs G/B)}
		\arrow["\circ", from=1-1, to=1-2]
		\arrow["\wr"', from=1-1, to=2-1]
		\arrow["\wr", from=1-2, to=2-2]
		\arrow["{(-,-)}", from=2-1, to=2-2]
	\end{tikzcd}\]
	where the pairing $(-,-)$ is $G_0(\pt/G)$-linear and admits an explicit description in terms of the Weyl character formula.
	In the other direction, composition is simply described via the box product
	
	\[\begin{tikzcd}
		{\Hom(p_*\un,\un)\times \Hom(\un,p_*\un)} & {\End(p_*\un)} \\
		{G_0(G\bs G/B)\times G_0(G\bs G/B)} & {G_0(G\bs(G/B\times G/B)).}
		\arrow["\circ", from=1-1, to=1-2]
		\arrow["\wr"', from=1-1, to=2-1]
		\arrow["\wr", from=1-2, to=2-2]
		\arrow["\boxtimes", from=2-1, to=2-2]
	\end{tikzcd}\]
	Moreover, the unit in $\End(p_*\un)$ corresponds to the class $i_*1\in G_0(G\bs(G/B\times G/B))$ for $1\in G_0(G\bs G/B)$ and $i:G\bs G/B\to G\bs(G/B\times G/B)$ the diagonal.

	Since $G$ has a simply-connected derived subgroup, the Pittie--Steinberg theorem \cite{steinbergTheoremPittie1975} shows that $\RR_G$ is a free $\RR_B$-module of rank $|W|$. Using this result, Kazhdan--Lusztig \cite[1.6 and 1.7]{kazhdanProofDeligneLanglandsConjecture1987} show that there are two bases $\{e_v\}_{v\in W}$ and $\{\hat{e}_v\}_{v\in W}$ of $\RR_B=G_0(G\bs G/B)$ as $\RR_G=G_0(G\bs \pt)$-modules such that $(e_v,\hat{e}_w)=\delta_{v,w}$. Moreover, they show that this yields a decomposition of the diagonal as $$i_*1=\sum_v e_v\boxtimes \hat{e}_v.$$
	In other words, the elements $e_v$ and $\hat{e}_v$ provide the projection and injection maps for the claimed direct sum decomposition $p_*\un=\un^{\oplus |W|}.$
\end{proof}
The results of Kazhdan--Lusztig \cite{kazhdanProofDeligneLanglandsConjecture1987} only consider the case of passing from $\pt/B$ to $\pt/G$. However, we will need a more general statement concerning parabolic subgroups. For now, we can show the following result, which will be improved in \Cref{cor:mapsofclassifyinspacesstronger}
\begin{corollary}\label{cor:mapsofclassifyinspaces}
	Let $G$ be a reductive group with simply-connected derived subgroup. Let $P\subset G$ be a parabolic subgroup. Denote by $p: \pt/P\to\pt/G$ the projection. Then $p_!(\un)=p_*(\un)$ is a direct summand of $\un^{\oplus n}\in \DK(\pt/G)$ for some $n\geq 1$.
\end{corollary}
\begin{proof}
	Denote by $B\subset P\subset G$ a Borel subgroup. Now consider the diagram
	
	\[\begin{tikzcd}
	{\pt/B} & {\pt/P} & {\pt/G.}
	\arrow["{p''}"', from=1-1, to=1-2]
	\arrow["{p'}", curve={height=-12pt}, from=1-1, to=1-3]
	\arrow["p"', from=1-2, to=1-3]
	\end{tikzcd}\]
	Using homotopy invariance, we can pass to the Levi factor of $P$ and apply \Cref{thm:pushforwardpittiesteinberg} to $p'$ and $p''$ to obtain that
	$$p_*\un^{\oplus|W_P|}=p_*p''_*\un=p'_*\un=\un^{\oplus|W_G|}.$$
	This shows that $p_*\un$ is a direct summand of $\un^{\oplus|W_G|}$.
\end{proof}
We obtain the following result concerning the functoriality.
\begin{lemma}\label{lem:FunctorialityOfDKT}
	Let $P\subset G$ be a parabolic subgroup of a reductive group with simply-connected derived subgroup.
	Then there are commutative diagrams
	\[\begin{tikzcd}
	{\DKT(\pt/P)} & {\op{Perf}(\RR_P)} \\
	{\DKT(\pt/G)} & {\op{Perf}(\RR_G)}
	\arrow[from=1-1, to=1-2]
	\arrow["{p^*\cong p^!}"', from=1-1, to=2-1]
	\arrow["{\otimes_{\RR_P} \RR_Q}", from=1-2, to=2-2]
	\arrow[from=2-1, to=2-2]
\end{tikzcd}
\qquad\text{and}\qquad
\begin{tikzcd}
	{\DKT(\pt/P)} & {\op{Perf}(\RR_P)} \\
	{\DKT(\pt/G)} & {\op{Perf}(\RR_G)}
	\arrow[from=1-1, to=1-2]
	\arrow["{p_*=p_!}"', from=1-1, to=2-1]
	\arrow["{\For }", from=1-2, to=2-2]
	\arrow[from=2-1, to=2-2]
\end{tikzcd}\]
where $p : \pt/P \to \pt/G$ is the natural map. 
\end{lemma}
\begin{proof}
	We first note that $p$ is smooth and proper. Hence $p_!=p_*$, and by Bott-periodicity $p^!\cong p^*$.
	First, we need to show that the functors $p_*$ and $p^*$ preserve the category $\DKT$. For the functor $p^*$ this is clear, since $p^*\un=\un$. For the functor $p_*$, this is shown in \Cref{cor:mapsofclassifyinspaces} where we use the assumption that the derived subgroup of $G$ is simply-connected.

	Now the second diagram commutes using the adjunction
	\begin{align*}
		\Map_{\DK(\pt/G)}(\un,p_*-)=\Map_{\DK(\pt/P)}(p^*\un,-)=\Map_{\DK(\pt/P)}(\un,-).
	\end{align*}
	Hence, the first diagram also commutes since both the vertical arrows are the left adjoints of the vertical arrows in the second diagram.
\end{proof}
\begin{remark}\label{remark:canonical-self-duality}
	\Cref{lem:FunctorialityOfDKT} implies that the induction and restriction functors $$ -\otimes_{\RR_P}\RR_G: \op{Perf}(\RR_P)\leftrightarrows \op{Perf}(\RR_G):\op{For}$$ are biadjoint. This reflects the fact that the inclusion $\RR_P\subset \RR_G$ is a so-called Frobenius extension.
\end{remark}
We can now show a stronger version of \Cref{cor:mapsofclassifyinspaces}.
\begin{corollary}\label{cor:mapsofclassifyinspacesstronger}
	Let $G$ be a reductive group with simply-connected derived subgroup. Let $P\subset G$ be a parabolic subgroup. Denote by $p: \pt/P\to\pt/G$ the projection. Then $p_!(\un)=p_*(\un)=\un^{\oplus |W/W_P|}$, where $W$ and $W_P$ denote the Weyl group of $G$ and $P$, respectively.
\end{corollary}
\begin{proof}
	Using \Cref{lem:FunctorialityOfDKT}, the object $\pi_* \un$ corresponds to the $\RR_G$-module $\RR_P$. Since the derived subgroup of $G$ is simply-connected, we can apply the Pittie--Steinberg theorem \cite{steinbergTheoremPittie1975} again to see that $\RR_P\cong {\RR_G}^{\oplus |W/W_P|}$. This implies the result.
\end{proof}
Taking these results, we are able to show the following generalization of the projective bundle formula.
\begin{corollary}\label{cor:projectivebundleformula}
	Under the assumptions of \Cref{cor:mapsofclassifyinspacesstronger}, assume that there is a Cartesian diagram of linearly reductive stacks
\[\begin{tikzcd}
	\Xx & \Yy \\
	{\pt/P} & {\pt/G}
	\arrow["\pi", from=1-1, to=1-2]
	\arrow["{q'}"', from=1-1, to=2-1]
	\arrow["q", from=1-2, to=2-2]
	\arrow["p", from=2-1, to=2-2]
\end{tikzcd}\]
	Then there is a natural equivalence of functors $\pi_*\pi^*\cong \un^{\oplus |W/W_P|}\otimes-$ on $\DK(\Yy)$.
\end{corollary}
\begin{proof}
	Using \Cref{cor:mapsofclassifyinspacesstronger} and base change, we obtain that $\pi_*\un=q^*p_*p^*\un=\un^{\oplus |W/W_P|}$.
	Now using the projection formula, we obtain natural equivalences of functors
	$$\pi_*\pi^*\cong \pi_*(\un\otimes \pi^*-)=\pi_*(\un)\otimes -=\un^{\oplus |W/W_P|}\otimes -$$
	which proves the claim.
\end{proof}

\begin{lemma}\label{lem:forget-along-diag}
Let $P$ be a parabolic subgroup of $G$ which is assumed to be of simply connected type. Consider the diagonal map
$$\Delta : \pt/P \to \pt/P \times \pt/P.$$
Then $\Delta_*$ restricts to a functor $\DKT(\pt/P)\to \DKT(\pt/P \times \pt/P)$ to Tate motives. 
\end{lemma} 

\begin{proof}
Consider the commutative diagram 
\[\begin{tikzcd}
	{\pt/B} & {\pt/P} \\
	{\pt/B \times \pt/B} & {\pt/P \times \pt/P.}
	\arrow[from=1-1, to=1-2]
	\arrow["\Delta"', from=1-1, to=2-1]
	\arrow["\Delta", from=1-2, to=2-2]
	\arrow[from=2-1, to=2-2]
\end{tikzcd}\]\\
By Theorem \ref{thm:pushforwardpittiesteinberg}, the pushforward along $\pt/B \to \pt/P$ (resp. $\pt/(B \times B) \to \pt/(P \times P)$) sends $\un$ to a direct sum of copies of $\un$. Hence it is enough to prove the statement for $P = B$. Consider now the diagram 
\[\begin{tikzcd}
	{\pt/T} & {\pt/B} \\
	{\pt/T \times \pt/T} & {\pt/B \times \pt/B.}
	\arrow["a", from=1-1, to=1-2]
	\arrow["\Delta"', from=1-1, to=2-1]
	\arrow["\Delta", from=1-2, to=2-2]
	\arrow["{a\times a}"', from=2-1, to=2-2]
\end{tikzcd}\]
Since $a_*\un = \un$ and $(a\times a)_*\un = \un$, it is enough to consider the case $B = T$. Using that $\pt/T = (T \times T)\backslash (T \times T)/T \to \pt/(T \times T)$, it is therefore enough to prove that $b_*\un$ is Tate, where $b : (T \times T)\backslash (T \times T)/T \to \pt/(T \times T)$. As $T$ is commutative, we can rewrite $(T \times T)/T = (\Gm)^r$ in a $(T \times T)$-equivariant way by choosing any basis of $X(T)$. The claim now follows from the definition of Tateness.
\end{proof}

\subsection{Non-representable functors}\label{sec:non-representable-functors}
Let $P,Q$ be parabolic subgroups of $G$ which is assumed to be of simply connected type.

\begin{construction}[Pushforward along a non-representable map]
We are going to need to study the pushforward along the map $\beta_P : \pt/P \to \pt$. As these maps are not defined in the $6$-functor formalism $\DK$, we only define it on Tate objects.
Namely, we consider the functor $\beta_*$ defined by
\[\begin{tikzcd}[column sep = 50pt]
	{\DKT(\pt/P)} & {\op{Perf}(R_P)} \\
	& {\Dd(\Z)}
	\arrow["{\Map(\un,-)}", from=1-1, to=1-2]
	\arrow["{\beta_*}"', from=1-1, to=2-2]
	\arrow[from=1-2, to=2-2]
\end{tikzcd}\]
\end{construction} 

\begin{lemma}\label{lem:Compatibility-Pullback}
Let $X$ be a stack equipped with two representable maps $X \xrightarrow{\alpha_1} \pt/(P \times Q)$ and $X \xrightarrow{\alpha_2} \pt/Q$ making the following diagram commute
\[\begin{tikzcd}
	X \\
	{\pt/(P \times Q)} & {\pt/Q}
	\arrow["{\alpha_1}"', from=1-1, to=2-1]
	\arrow["{\alpha_2}", from=1-1, to=2-2]
	\arrow["\pr"', from=2-1, to=2-2]
\end{tikzcd}\]
where the bottom map is simply the projection. Then there is a natural isomorphism of functors $\DKT(\pt/Q) \to \DK(X)$,
$$\alpha_1^*\pr^* \simeq \alpha_2^*.$$
\end{lemma}

\begin{proof}
Since $\DKT(\pt/Q) \simeq \op{Perf}(\RR_Q)$ and all categories in sight are $\RR_Q$-linear, it is enough to produce a morphism $\alpha_1^*\pr^*(\un) \simeq \alpha_2^*\un$. However we have $\alpha_1^*\pr^*(\un) = \un_X = \alpha_2^*\un$.
\end{proof}

\begin{remark}\label{remark:Compatibility-pushforward}
Passing to right adjoints in \Cref{lem:Compatibility-Pullback} yields an isomorphism of functors
$$\pr_*\Map(\un,\alpha_{1,*}-) \simeq \Map(\un, \alpha_{2,*}-).$$
\end{remark}

\section{Kac-Moody groups and flag varieties}\label{sec:Kac-Moody-groups}
We recall definitions and properties of Kac--Moody groups and their flag varieties. Moreover, we discuss a generalization to disconnected groups.
\subsection{Kac--Moody root data}\label{sec:kacmoodyrootdatum}
Let $\mathcal{D}=(X, \Delta, \Delta^\vee)$ be a Kac--Moody root datum with generalized Cartan matrix $A = (a_{\alpha,\beta})_{\alpha,\beta \in \Delta}$.
Here, $X$ is a finite-rank lattice, $\Delta\subset X$ and $\Delta^\vee\in X^\vee$ are finite subsets with a bijection $(-)^\vee:\Delta\to\Delta^\vee$, $a_{\alpha,\beta} = \langle \alpha, \beta ^{\vee}\rangle$, and $X^\vee=\Hom_\Z(X,\Z)$. 

A Kac--Moody root datum $\mathcal{D}$ is called 
\begin{enumerate}
	\item \emph{free} if the set of simple roots $\Delta$ is linearly independent, 
	\item \emph{of simply-connected type} if for all $\alpha\in \Delta$, $\langle -,\alpha^\vee\rangle: X\to \Z$ is surjective and
	\item \emph{of finite type} if $\det A\neq 0.$
\end{enumerate}
We denote the sets of (simple and positive) roots by $\Delta\subset \Phi_+\subset \Phi\subset X$.
\subsection{Kac--Moody groups and parabolic subgroups}\label{sec:kacmoodyrootdatum}
Let $\mathcal{D}$ be a Kac--Moody root datum. To this, we can associate the Kac--Moody group, Borel subgroup and maximal torus $G\supset B\supset T$. Here, $G$ is a group ind-scheme, $B$ is a pro-group and $T$ an algebraic group over $k.$ See, for example, \cite{mathieuConstructionGroupeKacMoody1989} and \cite{rousseauGroupesKacMoodyDeployes2016}.
We denote by $W=N(G)/T$ the Weyl group.

Each standard parabolic subgroup $G\supset P$ admits a Levi decomposition $P=L_P\ltimes U_P$ where $L_P$ is a Kac--Moody group corresponding to a Kac--Moody datum $\mathcal{D}_P$ and $U_P$ a pro-unipotent group. We denote the (simple and positive) roots, etc., associated to $\mathcal{D}_{P}$ by $\Delta_P, \Phi^+_P, \Phi_P,$ etc.
We say that $P$ is standard if it contains $B$, and that $P$ is of finite type if $\mathcal{D}_P$ is of finite type. 
In the following, \emph{we will only consider standard parabolic subgroups of finite type and will hence simply refer to them as parabolic subgroups.}

The unipotent radical $U_P\subset P$ of a parabolic $P$ is pro-unipotent, which can be seen as follows. For each $n\geq 0$, denote by $U_P(n)\subset U_P$ the subgroup generated by all root subgroups $U_\alpha\subset U_P$ of $P$-relative height $\geq n$; that is, when $\alpha$ is expressed in terms of simple roots $\Delta$, the sum of the coefficients of all $\Delta_P$ is $\geq n$. Then $U_P(n)$ is normal in $U_P$, the quotient $U_P/U_P(n)$ is finite dimensional, and the sequence of groups $U_P(n)$ equips $U_P$ with the structure of a pro-unipotent group. Moreover, the pro-unipotent structures of the unipotent radicals of parabolic subgroups $P\subset P'$ are compatible, see \cite[Section 6.1.13]{kumarKacMoodyGroupsTheir2002}. In particular, $U_P=U_{P'}\ltimes S$ for a finite-dimensional unipotent subgroup $S\subset P'$ normalized by $P$. 
\subsection{Weyl groups}

We denote by $W$ the Weyl group of $(G, T)$, and for parabolics $P,Q \subset G$, we denote by $W_P$ (resp. $W_Q$) the corresponding parabolic subgroup of $W$ attached to $P$ (resp. $Q$) and by ${^P}W$ (resp. $W^Q$, resp. ${^P}W^Q$) the set of minimal left $W_P$-coset representatives (resp. right $W_Q$-cosets, resp. $(W_P,W_Q)$-bicosets). In particular, we have ${^P}W^Q = {^P}W \cap W^Q$ and ${^P}W^B = {^P}W$, ${^B}W^Q = W^Q$. 

\subsection{Kac-Moody partial affine flag varieties}
Let $\mathcal{D}$ be a Kac--Moody datum and denote by $G\supset B\supset T$ the associated Kac--Moody group, Borel subgroup and maximal torus. Recall that in the following all parabolics are assumed to be standard and of finite type. In particular, quotients $Q/P$ of parabolics are finite-dimensional flag varieties.

We now briefly explain how double quotients $P\bs G/Q$ by parabolics are \niceindprostack ind-pro-stacks. The Bruhat decomposition shows that $G/Q$ is an ind-variety which arises as union of the $P$-stable varieties
$$(G/Q)_{\leq n}=\bigcup_{\substack{w\in W_P\bs W/W_Q\\ \ell(w)\leq n}} PwQ/Q.$$
Now, as in \cite[Section 2.3]{eberhardtUniversalKoszulDuality2025a}, for each $n\geq 0$ there is an $m(n)\geq 0$ such that $U_P(m(n))$ stabilizes $wQ/Q$ in $G/Q$ for all $w\in W$ with $\ell(w)\leq n$. Hence, the action of $P$ on $(G/Q)_{\leq n}$ factors over $P/U_P(m(n))$, so that $(P\bs G/Q)_{\leq n}$ becomes a \niceindprostack pro-stack where the transition maps are affine bundles. The maps $(P\bs G/Q)_{\leq n}\to (P\bs G/Q)_{\leq n+1}$ are pro-closed immersions. Hence $P\bs G/Q$ becomes a \niceindprostack ind-pro stack.

The Bruhat cells in $P\bs G/Q$ admit the following description in terms of classifying spaces.
\begin{lemma}\label{lem:P-orbits-are-nice}
For all $P,Q$ parabolics and $w \in \phantom{}^PW^Q$, the locally closed subscheme $P\bs PwQ/Q$ is an affine bundle over $\pt/Q_w$ where $Q_w$ is the parabolic associated to $W_P\cap wW_Qw^{-1}$.
\end{lemma}
\begin{proof}
Let $Q_w'=P\cap wQw^{-1}$ be the stabilizer of $wQ/Q$ in $P$. Then $Q_w'\subset Q_w$ and a standard result on intersections of (non-standard) parabolics says that $Q_w'R_u(Q_w)=Q_w.$
Hence, we obtain $P\bs PwQ/Q\cong \pt/Q_w'\to \pt/Q_w$, where the last map is an affine bundle.
\end{proof}

\subsection{Passing to disconnected groups}\label{sec:disconnected-groups}

We now consider the more general setup where we allow $G$ to be disconnected. Let $\mathcal{D}$ be a Kac-Moody root datum and let $G$ be a ind-pro-group scheme as above such that 
\begin{enumerate}
\item $G^{\circ}$ is the (connected) Kac-Moody group associated to $\mathcal{D}$, 
\item the group $\pi_0(G)$ is a discrete group ind-scheme.
\end{enumerate}

Examples of such disconnected groups are provided by loop groups $LH$ attached to finite-dimensional reductive (but maybe not semisimple) groups $H$. The group $H$ is not assumed to be connected either. Classical examples also include orthogonal groups $\mathrm{O}_n$ and normalizers of Levi subgroups. 
\begin{definition}
In a (possibly disconnected) group $G$, a parabolic subgroup $P \subset G$ is a parabolic subgroup of $G^{\circ}$; recall that these parabolics are always assumed to be of finite type.
\end{definition}

Let us fix $(B,T)$ a Borel pair of $G^{\circ}$. We introduce the following notations: 
\begin{enumerate}
\item we write $N^{\circ} = N_{G^{\circ}}(T)$ and $N = N_G(T)$,
\item $W^{\circ}$ is the Weyl group of $(G^{\circ},T)$, 
\item $W = N/T$ is the Weyl group of $G$,
\item $\Gamma = \pi_0(G)$ is the group of connected components, for a component $\gamma \in \Gamma$, we write $G^{\gamma}$ for the corresponding connected component of $G$. 
\end{enumerate}

\begin{lemma}\label{lem:splitting}
There is a short exact sequence 
$$1 \to W^{\circ} \to W \to \Gamma \to 1.$$
The choice of $B$ yields a splitting of this exact sequence. 
\end{lemma}

\begin{proof}
Since we have an exact sequence 
$$1 \to G^{\circ} \to G \to \Gamma \to 1$$ 
and $N^{\circ} = N \cap G^{\circ}$, we have an exact sequence 
$$1 \to N^{\circ} \to N \to N/N^{\circ} \to 1$$
and an inclusion $N/N^{\circ} \subset \Gamma$. It is enough to prove that $N \to \Gamma$ is surjective. Let $\gamma \in \Gamma$ and pick any $g \in G^{\gamma}$, then $^gT$ is a maximal torus of $G^{\circ}$ and there exists $h \in G^{\circ}$ such that $^{hg}T = T$. The element $hg \in G^{\gamma}$ lies in $N$ and maps to $\gamma$ in $\Gamma$. 

Let $\Gamma_B = N_G(T,B)/T$. It is clear that $\Gamma_B \cap W^{\circ} = 1$. We prove that $\Gamma_B$ surjects on $\Gamma$. Let $g \in G^{\gamma}$, and consider the pair $({^gT}, {^gB})$ then there exists $h \in G^{\circ}$ such that $({^{hg}T}, {^{hg}B}) = (T,B)$, i.e. $hg \in G^{\gamma}$ lies in $N_G(T,B)^{\gamma}$. 
\end{proof}

\begin{remark}\label{remark:multiplication-of-minimal-B-cosets}
The datum of $B$ yields a length function $\ell : W \to \N$ and the section of $W \to \Gamma$ is provided by the set of length $0$ elements. In particular since $B$ is normalized by $\gamma \in \Gamma$, the double coset $B\gamma B$ is isomorphic (as a left or right)-scheme to $B$ and we have, for all $\gamma_1, \gamma_2 \in \Gamma$, 
$$B\gamma_1 B \times^{B} B\gamma_2 B \xrightarrow{\sim} B \gamma_1\gamma_2 B.$$ 
\end{remark}

\begin{remark} \label{rem:comparisonofconnectedcomponents}
Let $P,Q$ be parabolics of $G^{\circ}$ and $\gamma \in \Gamma$. Via the splitting in \Cref{lem:splitting}, we may view $\gamma$ as an element in $W$ and we choose a representative $\dot\gamma\in  N$. Since $\gamma$ stabilizes $B$, $\gamma P \gamma^{-1}$ is also a standard parabolic.
Now left-multiplication and conjugation by $\dot\gamma$, respectively, yield isomorphisms $G^\circ\to G^\gamma$ and $P\to \gamma P \gamma^{-1}$ as well as an isomorphism
$$\dot{\gamma}:P\bs G^\circ/Q\to ( \gamma P \gamma^{-1})\bs G^\gamma /Q.$$
This is very useful to generalize statements from the connected to the disconnected case.
\end{remark} 

%% file: Inhalt/singularKtheoryBimodules.tex
\section{Singular $K$-theory Soergel bimodules}\label{sec:Singular-K-theory-SBim}
We define a combinatorial category of singular $K$-theory Soergel bimodules, which are a $K$-theoretic version and generalization to disconnected groups of singular Soergel bimodules, see \cite{williamsonSingularSoergelBimodules2010}. In this section, $\mathcal{D}$ is a free Kac--Moody root datum of simply-connected type. We denote by $G$ a (possibly disconnected) group with $G^{\circ}$ being the Kac--Moody group attached to $\mathcal{D}$. The group $G^{\circ}\supset B\supset T$ comes with a Borel subgroup and maximal torus. By convention, we consider only standard parabolic subgroups of finite type and simply refer to them as parabolics. Moreover, we denote $\Gamma = \pi_0(G)$.
\subsection{The Schur-algebroid}\label{sec:Schur-Algebroid}
In a first step, we define a category to encode the data of parabolics and connected components.
Namely, let $\Parab(G)$ be the category defined by
\begin{enumerate}
\item objects are parabolics of $G^{\circ}$,
\item for $P,Q$ parabolics of $G^{\circ}$, we set $\Hom_{\Parab(G)}(P,Q) = \Gamma$, and composition is given by multiplication in $\Gamma$.
\end{enumerate}
\begin{remark}
The category $\Parab(G)$ is a groupoid; the subgroupoid $\Parab(G^{\circ})$ is contractible and isomorphic to a point, and $\Parab(G) = \Parab(G^{\circ}) \times B\Gamma$. Hence, while categorically trivial, this category will be useful to organize the various Hecke categories.
\end{remark} 
\begin{remark} 
The Schur algebroid of $G^{\circ}$, following \cite{williamsonSingularSoergelBimodules2010}, is the functor
$$\Parab(G^{\circ}) \to \Morita(\Ab)$$
into the Morita category on $\Ab$, see Section \ref{sec:Morita-category}, defined as 
\begin{align*}
P &\mapsto \mathrm{H}_P \\
(P \to Q) &\mapsto \mathrm{H}_{P,Q}.
\end{align*}
where $\mathrm{H}_P$ is the parabolic Hecke algebra attached to $P$, and $\mathrm{H}_{P,Q}$ is the corresponding $(\mathrm{H}_P, \mathrm{H}_Q)$-bimodule. These are all free $\Z[q^{\pm 1}]$-modules. Specializing $q \to 1$ yields the variant with group algebras
\begin{align*}
P &\mapsto \Z[W_P\bs W/W_P] \\
(P \to Q) &\mapsto \Z[W_P \bs W/W_Q].
\end{align*} 
\end{remark}  

\subsection{The Bott--Samelson category}\label{subsec:Bott--Samelson-Category}
To organize the combinatorial data of Bott--Samelson resolutions of Schubert varieties in all $P\bs G^\gamma/Q$ simultaneously, we now define the \emph{Bott--Samelson category}, $\BS$, in the following way.
\begin{enumerate}
\item Objects are parabolic subgroups of $G^{\circ}$,
\item for $P,P'$ parabolic subgroups, we set $\Hom_{\BS}(P,P')$ to be the set of all sequences of parabolic subgroups
$$(P = P(0) \subset Q(0) \supset P(1) \subset Q(1) \subset \dots \subset Q(n) \supset P({n+1}) = P')$$
decorated by tuples $(\gamma_0, \dots, \gamma_n)$,
\item composition is given by concatenation of sequences.
\end{enumerate}
\begin{remark}
We will eventually write a sequence as above as 
$$(P = P(0) \subset (Q(0), \gamma_0) \supset P(1) \subset (Q(1), \gamma_1) \subset \dots \subset (Q(n), \gamma_n) \subset P({n+1}) = P').$$
If we want to shorten the notation, we will also write it as $(\underline{Q}, \underline{\gamma})$ or simply $(\underline{Q})$. 
\end{remark} 
The category $\BS$ is naturally fibered over $\Parab(G)$ via the functor 
\begin{align*}
\BS &\to \Parab(G) \\
P &\mapsto P \\
(\underline{Q}, \underline{\gamma}) &\mapsto \prod(\underline{\gamma}),
\end{align*}
where $\prod(\underline{\gamma})$ denotes the product $\gamma_0\gamma_1\dots \gamma_n$ in $\Gamma$. 

\subsection{Singular Soergel bimodules}
For a parabolic $P$, we denote the representation ring of $P$ by
$\RR_P = K_0(\pt/P)$. By assumption, $G$ is of simply-connected type, so that the derived subgroup of the Levi subgroup $L$ of $P$ is a simply-connected semisimple algebraic group. Recall that by the Chevalley restriction theorem, we have
$\RR_P = K_0(\pt/P)=\RR^W_P$
where $\RR=\RR_B=\Z[X(T)]$.
\subsubsection{Bott--Samelson and Soergel bimodules}\label{subsec:Bott--Samelson-And-Singular-Bimodules}
Let us proceed to define singular Bott--Samelson bimodules. Given a length $1$ sequence $(P, (Q, \gamma), P')$ in $\BS$, we write $Q^{\gamma}$ for the unique closed double $Q$-coset in $G^{\gamma}$, and we denote by $$\BS(P, (Q, \gamma), P') = K_0(P \backslash Q^{\gamma}/P')$$
\begin{lemma}\label{lem:Bott--Samelson-bimod-one-step}
Recall that, by assumption, $G$ is of simply-connected type. As an $\RR_P \otimes \RR_{P'}$-bimodule,
$$K_0(P \backslash Q/P') = \RR_{P} \otimes_{\RR_{Q}} \RR_{P'}.$$
\end{lemma} 

\begin{remark}
Lemma \ref{lem:Bott--Samelson-bimod-one-step} provides an explicit description of Bott--Samelson bimodules in the case $\gamma = 1$. 
\end{remark}

\begin{proof}
This can be computed geometrically in the following way. Consider the Cartesian square 
\[\begin{tikzcd}
	{\pt/P} & {P\backslash Q/P'} \\
	{\pt/Q} & {\pt/P'}
	\arrow["a"', from=1-1, to=2-1]
	\arrow["b"', from=1-2, to=1-1]
	\arrow["{b'}", from=1-2, to=2-2]
	\arrow["{a'}", from=2-2, to=2-1]
\end{tikzcd}\]\\
By Lemma \ref{lem:DescriptionOfDKT}, we have $K_0(P\backslash Q/P') = \End_{\DK(P\backslash Q/P')}(\un)$ and using that $\un_{P\backslash Q/P'} = b^*\un_{\pt/P} = b'^{*}\un$, we have by adjunction and proper base change 
\begin{align*}
\End_{\DK(P\backslash Q/P')}(\un) &= \Hom_{P\backslash Q/P'}(b^*\un, b'^{*}\un) \\
&= \Hom_{\pt /P'}(\un_{\pt/P'},a'^*a_*\un_{\pt/P}).
\end{align*}
The result now follows from Lemma \ref{lem:FunctorialityOfDKT}.
\end{proof}

\begin{lemma}
Let $\gamma \in \Gamma$. There is a natural isomorphism as a left (resp. right) $\RR$-module
$$K_0(B \backslash B\gamma B/B) \simeq \RR$$
\end{lemma}

\begin{proof}
This follows from the fact that $B \backslash B\gamma B/B \simeq \pt/B$ by either using the left or right projection to $\pt/B$, see Remark \ref{remark:multiplication-of-minimal-B-cosets}. 
\end{proof}
\begin{remark}
We do not give a formula for $K_0(P \backslash P\gamma Q/Q)$ in general. 
\end{remark} 

 This is naturally a $(\RR_P, \RR_{P'})$-bimodule. Consider now the following functor 
\begin{align*}
\BS &\to \Morita(\Ab) \\
P &\mapsto \RR_P \\
(\underline{Q}, \underline{\gamma}) &\mapsto \BSBim(\underline{Q}, \underline{\gamma}),
\end{align*}
where $\BSBim(\underline{Q}, \underline{\gamma})$ is defined as the tensor product
\begin{align*}
    \BS(P(0), (Q(0), \gamma_0), P(1)) \otimes_{\RR_{P(1)}} \dots \otimes_{\RR_{P(n)}} \BS(P(n), (Q(n), \gamma_n), P(n+1)).
\end{align*}
In particular, if $\gamma_i=1$, we obtain that
$$\BSBim(\underline{Q})=\RR_{P(0)}\otimes_{\RR_{Q(0)}}\RR_{P(1)}\otimes_{\RR_{Q(1)}}\dots \otimes_{\RR_{Q(n)}}\RR_{P(n+1)}.$$
\begin{remark}
We invite the reader to compare with the usual construction of Bott--Samelson bimodules, see \cite{soergelKategorieMathcalPerverse1990,williamsonSingularSoergelBimodules2010,soergelEquivariantMotivesGeometric2018}. A difference with our setup is that the ring $\RR$ is not graded and that we allow for disconnected groups.
\end{remark}

Given $(P,Q, \gamma)$, we denote by $\BSBim_{P,Q}^{\gamma} \subset \Mod_{\RR_P \otimes_{\Z} \RR_Q}$ the full additive subcategory whose objects are $\BSBim(\underline{Q}, \underline{\gamma})$ for all sequences from $P$ to $Q$ such that $\prod \underline{\gamma} = \gamma$.

\begin{remark}
All the categories $\BSBim_{P,Q}^{\gamma}$ are $1$-categories. 
\end{remark} 

\begin{definition}
The category $\SBim_{P,Q}^{\gamma}$ is the Karoubian completion of $\BSBim_{P,Q}^{\gamma}$. We can sum together the categories and set 
$$\SBim_{P,Q} = \bigoplus_{\gamma \in \Gamma} \SBim_{P,Q}^{\gamma}.$$
This is a full subcategory of $\Gamma$-graded $\RR_P \otimes \RR_Q$-modules. 
\end{definition} 

\begin{remark}
It follows from the construction that $\SBim_{P,P}$ and $\SBim_{P,P}^{\circ}$ are naturally monoidal categories (for the convolution of bimodules) and $\SBim_{P,Q}$ (resp. $\SBim_{P,Q}^{\gamma}$) is naturally a $\SBim_{P,P}-\SBim_{Q,Q}$-bimodule (resp. a $\SBim_{P,P}^{\circ}-\SBim_{Q,Q}^{\circ}$-bimodule).
\end{remark} 

\begin{remark}
Note that none of the rings $\RR, \RR_P$ or $\RR_Q$ are local; the Krull-Schmidt theorem fails, and we cannot \emph{a priori} say that any of these objects are finite sums of indecomposable bimodules in a unique way.
\end{remark}

\begin{remark}
The functor 
\begin{align*}
\Parab(G) &\to \Morita(\mathrm{Pr}^{\omega}) \\
P &\mapsto \Ch^b(\SBim_{P,P}^{\circ})\\
P \xrightarrow{\gamma} Q &\mapsto \Ch^b(\SBim_{P,Q}^{\gamma}), 
\end{align*}
encapsulates the structure of singular $K$-theory Soergel bimodules. In \cite{williamsonSingularSoergelBimodules2010}, the theory is rather formulated in terms of monoidal $2$-categories, the two points of view are exchanged by the Grothendieck construction as pointed out in Remark  \ref{rmk:Grothendieck-construction}.
\end{remark}

%% file: Inhalt/ktheoreticsoergel2.tex
\section{Soergel description of the $K$-motivic Hecke category}\label{sec:Main-section}

In this section, we fix $G$ a Kac--Moody group of simply connected type with Borel pair $(B,T)$, we shall use the notations of Sections \ref{sec:Kac-Moody-groups} and \ref{sec:Singular-K-theory-SBim}. All parabolics $P,Q,\dots \subset G^{\circ}$ are assumed to be standard and of finite type. 

\subsection{$K$-motivic Hecke categories}

Recall that $k$ is an algebraically closed field of characteristic $0$. For parabolics $P,Q \subset G$, the \emph{Hecke stack} (of type $(P,Q)$) is $P \backslash G/Q$. It is a linearly reductive ind-pro stack and the category $\DK(\PGQ)$ is well defined. The stack $\PGQ$ admits a Whitney stratification 
$$\PGQ = \bigsqcup_{w \in {^PW^Q}} P\backslash PwQ/Q$$
and each locally closed substack $P\backslash PwQ/Q$ is a (pro)-affine bundle over a finite dimensional partial flag variety, see Lemma \ref{lem:P-orbits-are-nice}. 
For $w \in P\backslash PwQ/Q$, we shall denote the inclusion by
$$i_w^{P,Q} : P\backslash PwQ/Q \hookrightarrow \PGQ.$$
We will drop the exponent $^{P,Q}$ if it is clear from context. 
The standard and costandard objects are defined as
\begin{enumerate}
\item $\Delta_w^{P,Q} = i_{w, !}^{P,Q} \un$ and
\item $\nabla_{w}^{P,Q} = i_{w,*}^{P,Q} \un$. 
\end{enumerate}

\begin{definition}\label{def:heckecat}
The $K$-theoretic singular Hecke category of type $(P,Q)$ is the stable, thick subcategory generated by standard objects
$$\Hecke_{P,Q} = \langle \Delta_w^{P,Q}, w \in {^P}W^Q \rangle_{\op{stable},\op{thick}} \subset \DK(\PGQ).$$ 
If we want to distinguish the connected components, we will write 
$$\Hecke_{P,Q}^{\gamma} \subset \Hecke_{P,Q}$$
for the full subcategory of $K$-motives supported on $P\backslash G^{\gamma}/Q$.
\end{definition} 

\begin{remark}
\begin{enumerate}
\item As explained in \cite[Remark 4.2]{eberhardtUniversalKoszulDuality2025a}, this category is much smaller than $\DK(\PGQ)$.
\item Clearly, we have $\Hecke_{P,Q} = \bigoplus_{\gamma} \Hecke_{P,Q}^{\gamma}$.
\end{enumerate}
\end{remark} 
\subsection{Change of Parabolics}\label{sec:changeofparabolics}
For parabolics $P \subset P'$ and $Q \subset Q'$ we denote the projection map 
$$^{P}_{P'}\pi_{Q'}^{Q} : P \backslash G/Q \to P' \backslash G/Q'.$$
If $P=P'$, we simply write $\pi_{Q'}^{Q}$ and if $Q = Q'$, we write $^{P}_{P'}\pi$. This way, $^{P}_{P'}\pi_{Q'}^{Q}=\phantom{}^{P}_{P'}\pi\circ \pi_{Q'}^{Q}$.
We collect some important properties of these maps. We will often abbreviate $\pi=\phantom{}^{P}_{P'}\pi_{Q'}^{Q}$ in the following.
\begin{remark}
	The map $\pi$ restricts to connected components $P\bs G^\gamma/Q\to P'\bs G^\gamma/Q'$ for $\gamma\in \Gamma$. Moreover, the description $P\bs G^\circ /Q\cong (\gamma P\gamma^{-1})\bs G^\gamma/Q$ of any connected component in terms of the identity component is compatible with $\pi$. In the following results we can hence always safely assume that $G$ is connected.
\end{remark}
\begin{lemma}\label{lem:basicpropertyofchangeofparabolic}
The map $\pi$  is a locally trivial $P'/P \times Q'/Q$-fibration so that $\pi_*\cong \pi_!$ and $\pi^*\cong\pi^!.$	
\end{lemma}
Our first goal is to show that pushforward along $\pi$ preserves (up to direct sums and summands) standard and costandard objects. For this, we will use the following result.
\begin{lemma}\label{lem:ChangeOfParabolicsPushforwardOnBruhatCell}
	Let $w\in W^Q$ and denote by $w'\in \phantom{}^PW^Q$ the shortest representative of $W_PwW_Q$. Denote by $p: B\bs BwQ/Q\to P\bs Pw'Q/Q$ the projection. Then $p_*(\un)\cong p_!(\un)= \un^{\oplus n}$ for some $n\geq 1$.
\end{lemma}
\begin{proof}
	Write $xw'=w$, so that $\ell(w)=\ell(x)+\ell(w')$. By \cite[Theorem 3.2.1]{eberhardtGradedGeometricParabolic2018}, since $w'$ is a shortest coset representative the subgroup $W_P\cap w'W_Qw'^{-1}\subset W$ is generated by simple reflections. The associated parabolic subgroup $Q_{w'}\subset G$ acts on $Bw'Q$ by left multiplication and there is a commutative diagram
	
\[\begin{tikzcd}[column sep=10pt]
	{\pt/xBx^{-1}\cap Q_{w'}} & {B\bs BxQ_{w'}/Q_{w'}} & {B\bs BxQ_{w'}\times^{Q_{w'}}Bw'Q/Q} & {B\bs BwQ/Q} \\
	{\pt/Q_{w'}} & {P\bs P/Q_{w'}} & {P\bs P\times^{Q_{w'}}Bw'Q/Q} & {P\bs Pw'Q/Q}
	\arrow["\sim", from=1-1, to=1-2]
	\arrow["p'", from=1-1, to=2-1]
	\arrow[from=1-2, to=2-2]
	\arrow[from=1-3, to=1-2]
	\arrow["\sim", from=1-3, to=1-4]
	\arrow[from=1-3, to=2-3]
	\arrow[from=1-4, to=2-4]
	\arrow["\sim", from=2-1, to=2-2]
	\arrow[from=2-3, to=2-2]
	\arrow["\sim", from=2-3, to=2-4]
\end{tikzcd}\]
where the middle square is Cartesian. Now $Q_{w'}'=(xBx^{-1}\cap Q_{w'})R_u(Q_{w'})\subset Q_{w'}$. Hence, the map $p'$ factors as an affine bundle $\pt/(xBx^{-1}\cap Q_{w'})\to \pt/Q_{w'}'\to \pt/Q_{w'}$. The result follows from \Cref{cor:mapsofclassifyinspacesstronger} using that the Levi factors of $Q_{w'}'$ and $Q_{w'}$ have simply-connected derived subgroups by assumption.
\end{proof}
\begin{remark}
	In the case that $Q=B$ the map $p:B\bs BwB/B\to P\bs Pw'B/B$ is simply an affine bundle and $p_*(\un)=p_!(\un)=\un$.
\end{remark}
\begin{lemma}\label{lem:ChangeOfParabolicPushforward}
	Let $w\in \phantom{}^{P}W^{Q}$ and denote by $w'\in \phantom{}^{P'}W^{Q'}$ the shortest representative of $W_{P'}wW_{Q'}$. Then
	$\pi_*\Delta_w^{P,Q}$ is a direct summand of $(\Delta_{w'}^{P',Q'})^{\oplus n}$ and $\pi_*\nabla_w^{P,Q}$ is a direct summand of $(\nabla_{w'}^{P',Q'})^{\oplus n}$ for some $n\geq 1.$
\end{lemma}
\begin{proof}
	Since $\pi=\phantom{}^{P}_{P'}\pi_{Q'}^{Q}$ and the statement is symmetric in $P$ and $Q$, it is enough to treat the case $P=P'$. Abbreviate $\pi'=\pi_{Q,*}^B$, $\pi''=\pi_{Q',*}^B$. By \Cref{lem:ChangeOfParabolicsPushforwardOnBruhatCell}, we obtain $(\Delta_w^{P,Q})^{\oplus n}=\pi_{*}'\Delta^{P,B}_w$ for some $n\geq 1$. By using \Cref{lem:ChangeOfParabolicsPushforwardOnBruhatCell} again, we obtain
	$$\pi_{*}(\Delta_w^{P,Q})^{\oplus n}=\pi_{*}\pi'_{*}\Delta_{w}^{P,B} = \pi_{*}''\Delta^{P,B}_w = (\Delta_{w'}^{P,Q'})^{\oplus m}$$
	for some $m\geq 0$ so that $\pi_{*}(\Delta_{w}^{P,Q})$ is a direct summand of $(\Delta_{w'}^{P,Q'})^{\oplus m}$.
	The statement for costandard objects is dual.
\end{proof}
Next, we study the effect of pullback along $\pi.$
\begin{lemma}\label{lem:ChangeOfParabolicPullback}
	Let $w'\in \phantom{}^{Q'}W^{P'}.$
	Then $\pi^*\Delta_{w'}^{P',Q'}$ admits a cofiltration whose associated graded are objects $\Delta_w^{P,Q}$ for
	$w \in {^{P}}W^{Q}\cap W_{P'}w'W_{Q'}$.
\end{lemma}
\begin{proof}
	By base change $\pi^{*}\Delta_{w'}^{P',Q'}=j_{w',!}\un$ where we denote by
$$j_{w'} : P \bs P' w' Q' /Q \hookrightarrow P \bs G/Q.$$ Now $ P \bs P' w' Q'/Q$ admits a stratification with strata of the form $PwQ$ for $w \in {^{P}}W^{Q}\cap W_{P'}w'W_{Q'}$. The desired cofiltration arises by an iterated application of the localization triangle. 
\end{proof}
\begin{corollary} \label{cor:ChangeOfParabolic}
	The functors $\pi^{*}$ and $\pi_{*}$ induce a pair of bi-adjoint functors $$\pi_{*}:\Hecke_{P,Q} \leftrightarrows \Hecke_{P',Q'}:\pi^{*}.$$
\end{corollary}

\begin{corollary}\label{corol:Generation-By-BB}
The category $\Hecke_{P,Q}$ is stably generated by direct summands of objects of the form $\pi_*A$ for $A \in \Hecke_{B,B}$ where $\pi={^B_P}\pi^{B}_{Q}$.
\end{corollary}
\subsection{Monoidal structures}

We introduce the usual monoidal structures. Let $P_1,P_2$ and $P_3$ be three parabolics and consider the diagram 
\[\begin{tikzcd}[column sep=0]
	& {P_1 \backslash G \times^{P_2} G/P_3 } && {P_1 \backslash G/P_3} \\
	{P_1 \backslash G/P_2 } && {P_2 \backslash G/P_3}
	\arrow["m", from=1-2, to=1-4]
	\arrow["{\pr_1}", from=1-2, to=2-1]
	\arrow["{\pr_2}"', from=1-2, to=2-3]
\end{tikzcd}\]
We have two bifunctors 
$$\DK(P_1 \backslash G/P_2) \times \DK(P_2 \backslash G/P_3) \to \DK(P_1 \backslash G/P_3)$$
given by 
$$A *^{P_2}_! B = m_!(\pr_1^*A \otimes \pr_2^*B)$$
and 
$$A *^{P_2}_* B = m_*(\pr_1^!A \otimes \pr_2^!B).$$
Note that the map $\pr_1$ (resp. $\pr_2$, resp. $m$) is ind-representable with fibers isomorphic to $G/P_3$ (resp. $G/P_1$, resp. $G/P_2$). This equips the category $\DK(P\backslash G/P)$ with two monoidal structures and the category $\DK(P\backslash G/Q)$ with two bimodule structures over $\DK(P\backslash G/P)$ and $\DK(Q \backslash G/Q)$.
\begin{lemma}\label{lem:comparison-convol-finite-dim}
Assume that $G$ is finite dimensional, then $*^{P}_! = *_*^P$. 
\end{lemma}
\begin{proof}
Since $G$ is finite dimensional, the map $m$ is proper, hence $m_! = m_*$. Since $G$ is finite dimensional, $G/P$ is a finite dimensional smooth scheme and therefore $\pr_i^* = \pr_i^!$ (see Remark \ref{remark:comparison-pull-back}).
\end{proof}

\begin{remark}
We only state Lemma \ref{lem:comparison-convol-finite-dim} in the finite dimensional case for now. The general case is proven in Theorem \ref{thm:Rigitidy-Of-Hecke-Categorioid}. As $G/P_2$ is ind-proper, we have $m_! = m_*$, the problem in comparing the two definitions geometrically comes from the $\pr_i^!$ vs. $\pr_i^*$ in the definition of the two convolutions, as $G/P_2$ is not ind-smooth, we cannot immediately compare them. The same issue arose in \cite{eberhardtUniversalKoszulDuality2025a}. 
\end{remark} 

\begin{remark}
Since the parabolics $P_1,P_2$ and $P_3$ are assumed to be contained in $G^{\circ}$, both convolutions define bifunctors 
$$\DK(P_1 \backslash G^{\gamma_1} /P_2) \times \DK(P_2 \backslash G^{\gamma_2}/P_3) \to \DK(P_1 \backslash G^{\gamma_1\gamma_2}/P_3).$$
\end{remark}

\begin{lemma}\label{lem:Action-Of-Min-Elts}
For any $\gamma \in \Gamma$ and $w \in {^P}W^B$ (resp. $w \in {^B}W^Q$), there are canonical isomorphisms 
\begin{enumerate}
\item $\Delta_w^{P,B} *_!^B \Delta_{\gamma}^{B,B} = \Delta_{w\gamma}^{P,B} = \Delta_w^{P,B} *_*^B \Delta_{\gamma}^{B,B}$,
\item $\Delta_{\gamma}^{B,B} *_!^B \Delta_w^{B,Q}  = \Delta_{\gamma w}^{B,Q} = \Delta_{\gamma}^{B,B} *^B \Delta_w^{B,Q}$.
\end{enumerate}
\end{lemma}

\begin{proof}
The two statements are symmetric, let us prove the second one. This is a direct consequence of the fact that 
$$B \backslash B^{\gamma} \times^B BwQ/Q \to B \backslash B\gamma wQ/Q$$ 
is an isomorphism (see Remark \ref{remark:multiplication-of-minimal-B-cosets}). 
\end{proof}

\subsection{Bott--Samelson objects}\label{sec:Bott-Samelson}
We now introduce the Bott--Samelson objects and discuss how convolution with Bott--Samelson objects is related to pullback and pushforward along the maps $^{P}_{P'}\pi_{Q'}^{Q} : P \backslash G/Q \to P' \backslash G/Q'$ discussed in \Cref{sec:changeofparabolics}. Our discussion is just a variation on \cite[Lemma 3.2.1]{soergelRelationIntersectionCohomology2000a}.
We will use the notations of Section \ref{subsec:Bott--Samelson-Category} for the Bott--Samelson category.
\begin{definition}
Let $P,P'$ be two parabolics of $G^{\circ}$ and let $(\underline{Q}, \underline{\gamma})$ be a map in $\BS$ from $P$ to $P'$. The Bott--Samelson map is 
$$\pi_{\underline{Q}} : P \backslash Q(0)^{\gamma_0} \times^{P(1)} Q(1)^{\gamma_1} \times^{P(2)} \dots \times^{P({n})} Q({n})^{\gamma_{n}} /P' \to P\backslash G/P'.$$
The Bott--Samelson-object $\BS(\underline{Q})$ is 
$$\BS(\underline{Q}) = \pi_{\underline{Q},*}\un.$$ 
\end{definition}

\begin{lemma}\label{lem:Basics-of-Bott--Samelson}
\begin{enumerate}
\item Let $(\underline{Q}, \underline{\gamma})$ be a sequence from $P$ to $P'$ and assume that $\gamma_i = 1$ for all $i$, then there is a canonical isomorphism of functors $$- *_!^{P} \BS(\underline{Q}) =  \pi_{Q(0)}^{P(1),*}\pi_{Q(0),*}^{P(1)} \dots \pi_{Q(n)}^{P',*}\pi_{Q(n),*}^{P'},$$ and similarly for $\BS(\underline{Q}) *_!^{P'} -$.
\item Let $\underline{Q}_1$ and $\underline{Q}_2$ be two sequences from $P_1 \to P_2$ and $P_2 \to P_3$ respectively, then for $? \in \{!,*\}$
$$\BS(\underline{Q}_1) *^{P_2}_? \BS(\underline{Q}_2) = \BS(\underline{Q}_1 \sqcup \underline{Q}_2).$$
\item All objects $\BS(\underline{Q})$ are objects of $\Hecke_{P,P'}$ (and not just in $\DK(P \backslash G/P')$), in particular convolving (using $*_!$ or $*_*$) with a Bott--Samelson preserves the category $\Hecke_{P,P'}$. 
\end{enumerate}
\end{lemma}

\begin{proof}
We start by proving item $(2)$. We will abbreviate $\underline{Q} = Q(0)^{\gamma_0} \times^{P(1)} \dots \times^{P(n)} Q(n)^{\gamma_n}$ and similarly for $\underline{Q}'$.
Consider the diagram 
\[\begin{tikzcd}[column sep=0pt]
	&& {P_1\backslash \underline{Q} \times^{P_2}\underline{Q}'/P_3} \\
	& {P_1\backslash \underline{Q} \times^{P_2}G/P_3} && {P_1\backslash G \times^{P_2} \underline{Q}'/P_3} \\
	{P_1\backslash \underline{Q}/P_2} && {P_1\backslash G \times^{P_2}G/P_3} && {P_2\backslash \underline{Q}'/P_3} \\
	& {P_1\backslash G/P_2} && {P_2\backslash G/P_3}
	\arrow["{\tilde{\pi}_{\underline{Q}'}}"', from=1-3, to=2-2]
	\arrow["{\tilde{\pi}_{\underline{Q}}}", from=1-3, to=2-4]
	\arrow["{\tilde{\pi}_{\underline{Q} \sqcup \underline{Q}'}}", from=1-3, to=3-3]
	\arrow["{\pr_1}"', from=2-2, to=3-1]
	\arrow["{\pi_{\underline{Q}}}"', from=2-2, to=3-3]
	\arrow["{\pi_{\underline{Q}'}}", from=2-4, to=3-3]
	\arrow["{\pr_2}", from=2-4, to=3-5]
	\arrow["{\pi_{\underline{Q}}}"', from=3-1, to=4-2]
	\arrow["{\pr_1}", from=3-3, to=4-2]
	\arrow["{\pr_2}"', from=3-3, to=4-4]
	\arrow["{\pi_{\underline{Q}'}}", from=3-5, to=4-4]
\end{tikzcd}\]
All three squares are Cartesian and by proper base change we have 
$$\pr_1^*\BS(\underline{Q}) \otimes \pr_2^*\BS(\underline{Q}') = \tilde{\pi}_{\underline{Q}\sqcup \underline{Q}',*}\un.$$
Hence we have $\BS(\underline{Q}) *_!^{P_2} \BS(\underline{Q}') = m_* \tilde{\pi}_{\underline{Q}\sqcup \underline{Q}',*}\un = \BS(\underline{Q} \sqcup \underline{Q}')$ since $m \circ \tilde{\pi}_{\underline{Q}\sqcup \underline{Q}'} = \pi_{\underline{Q} \sqcup \underline{Q}'}$. Since the maps $\pi_{\underline{Q}}$ and $\pi_{\underline{Q}'}$ are proper, we also have 
$$\BS(\underline{Q}) *_*^{P_2} \BS(\underline{Q}') = \BS(\underline{Q}\sqcup \underline{Q}').$$

We now prove item $(1)$. Using item $(2)$, we can proceed inductively and assume that the sequence $\underline{Q}$ has length one. Let us assume that $\underline{Q} = (P \subset (Q, 1) \supset P')$. Let $P_0$ be some parabolic and consider the diagram 
\[\begin{tikzcd}
	&& {P_0 \backslash G/P'} \\
	{P_0 \backslash G/P} & {P_0 \backslash G \times^{P} G/P'} & {P_0 \backslash G \times^{P} Q/P'} \\
	& {P \backslash G /P'} & {P \backslash Q /P'.}
	\arrow["m", from=2-2, to=1-3]
	\arrow["{\pr_1}"', from=2-2, to=2-1]
	\arrow["{\pr_2}"', from=2-2, to=3-2]
	\arrow["m"', from=2-3, to=1-3]
	\arrow["{\pi_Q}"', from=2-3, to=2-2]
	\arrow["{\pr_2}", from=2-3, to=3-3]
	\arrow["{\pi_Q}", from=3-3, to=3-2]
\end{tikzcd}\]
Since the bottom square is Cartesian, we have an isomorphism of functors 
$$- *^P_! \BS(Q) = m_*\pi_Q^*\pr_1^*(-).$$
Applying base change to the right hand side using the following Cartesian diagram shows item $(1)$.
\[\begin{tikzcd}[column sep = 0]
	& {P_0 \backslash G \times^{P}Q/P'} \\
	{P_0 \backslash G/P} && {P_0 \backslash G/P'} \\
	& {P_0 \backslash G/Q}
	\arrow["{\pr_1}"', from=1-2, to=2-1]
	\arrow["m", from=1-2, to=2-3]
	\arrow["{\pi_{Q}^P}"', from=2-1, to=3-2]
	\arrow["{\pi_{Q}^{P'}}", from=2-3, to=3-2]
\end{tikzcd}\]\\
Finally for item $(3)$, it is enough to consider the case of a length $1$ sequence, so we assume that $\underline{Q} = (P \subset (Q, \gamma) \supset P')$. If $\gamma = 1$, this is clear by item $(1)$ and Lemmas \ref{lem:ChangeOfParabolicPushforward} and  \ref{lem:ChangeOfParabolicPullback}. In general, the object $\BS(P \subset (Q,\gamma) \supset P')$ is a direct summand of 
$$\BS(P \subset (Q,1) \supset B) *^{B}_! \BS(B \subset (B,\gamma) \supset B) *_!^B \BS(B \subset (P',1) \supset P')$$
by Remark \ref{remark:Example-Bott--Samelson}, Lemma \ref{lem:Action-Of-Min-Elts} and Lemma \ref{lem:ChangeOfParabolicPushforward}.
\end{proof}
\begin{remark}\label{remark:Example-Bott--Samelson}
Let us highlight some special cases of Bott--Samelson objects.
\begin{enumerate}
\item Consider the sequence $P \subset (P, \gamma) \supset P$ for some $\gamma \in \Gamma$, then the corresponding Bott--Samelson object is $\Delta_{\gamma}^{P,P}$.
\item Consider the sequence $P \subset (P, 1) \supset P'$, for $P' \subset P$, then the corresponding Bott--Samelson object is $\Delta_1^{P,P'}$ and $\Delta_1^{P,P'} *^{P'} -$ is isomorphic to ${^{P'}_P}\pi_*$ and $- *^P \Delta_1^{P,P'}$ is isomorphic to $\pi^{P',*}_P$. 
\item Dually to case $(2)$, attached to the sequence $P' \subset (P,1) \supset P$, the corresponding bimodule is $\Delta_1^{P',P}$ and convolution on the left (resp. right) with $\Delta_1^{P',P}$ is isomorphic to $\pi^{P'}_{P,*}$ (resp. ${^{P'}_P}\pi^*$). 
\end{enumerate}
\end{remark} 

\begin{lemma}\label{lem:Generateting-families}
The category $\Hecke_{P,Q}$ is the thick stable subcategory generated by either three collections of objects:
\begin{enumerate}
\item the standard objects $(\Delta_w^{P,Q})_{w \in {^P}W^Q}$, 
\item the costandard objects $(\nabla_w^{P,Q})_{w \in {^P}W^Q}$,
\item the Bott--Samelson objects $\BS(\underline{Q})$ for all sequences $\underline{Q}$ from $P$ to $Q$. 
\end{enumerate}
\end{lemma} 

\begin{proof}
The first family generates $\Hecke_{P,Q}$ by definition. By \cite[Corollary 4.9]{eberhardtUniversalKoszulDuality2025a}, the category $\Hecke_{B,B}^{\circ}$ is generated by Bott--Samelson objects attached to sequences 
$$(B \subset P_{s_1} \supset B \subset P_{s_2} \supset \dots \subset P_{s_n} \supset B)$$
where $(s_1, \dots, s_n)$ is a sequence of simple reflections and $P_{s_i}$ is the minimal parabolic attached to $s_i$. Since all objects of the form $\Delta_{\gamma}^{B,B}$ are Bott--Samelson objects by Remark \ref{remark:Example-Bott--Samelson}, we see that the category $\Hecke_{B,B}$ is generated by Bott--Samelson objects. By Corollary \ref{corol:Generation-By-BB}, the category $\Hecke_{P,Q}$ is generated by pushforwards of objects in $\Hecke_{B,B}$, and by Remark \ref{remark:Example-Bott--Samelson}, this pushforward can be realized as convolution against some Bott--Samelson. It follows that the collection of all Bott--Samelson attached to sequences from $P$ to $Q$ generate $\Hecke_{P,Q}$. 

Given a sequence $\underline{Q}$ from $P$ to $Q$, since the map $\pi_{\underline{Q}}$ is proper and the stack
$$P \backslash Q(0)^{\gamma_0} \times^{P(1)} Q(1)^{\gamma_1} \times^{P(2)} \dots \times^{P({n})} Q(n)^{\gamma_{n}} /P' \to P\backslash G/P'$$
is (pro)-smooth, we have $\Dual(\BS(\underline{Q})) = \BS(\underline{Q})$. Since these objects generate $\Hecke_{P,Q}$, it follows that $\Hecke_{P,Q}$ is stable under Verdier duality. Finally, $\Dual(\Delta_w^{P,Q}) = \nabla_{w}^{P,Q}$, so the collection of costandard objects generate $\Hecke_{P,Q}$.
\end{proof}

\begin{lemma} 
The convolution $*_!^P$ and $*_*^P$ send $\Hecke \times \Hecke$ to $\Hecke$. 
\end{lemma} 

\begin{proof}
By Lemma \ref{lem:Generateting-families}, it is enough to prove that for two Bott--Samelson objects, their convolution (using $*_!$ or $*_*$) lies in $\Hecke$. This is clear by Lemma \ref{lem:Basics-of-Bott--Samelson}.
\end{proof}

\begin{construction}
Recall the definition of the category $\Parab(G)$, from Section \ref{sec:Schur-Algebroid}. We define a functor to organize the various Hecke categories as follows 
\begin{align*}
\Hecke : \Parab(G) &\to \Morita(\mathrm{Pr}^{\omega}) \\
P &\mapsto \Hecke_{P,P}^{\circ} \\
(P \xrightarrow{\gamma} Q) &\mapsto \Hecke_{P,Q}^{\gamma},
\end{align*}
where the composition is induced by $*_!$-convolution. 
\end{construction}

\begin{definition}
The total Hecke category $\Hecke_{\tot}$ is
$$\Hecke_{\tot} = \bigoplus_{P,Q} \Hecke_{P,Q}.$$
We view it as a monoidal category using $*_!$ and declaring that for $A \in \Hecke_{P_1,P_2}$ and $B \in \Hecke_{P_3,P_4}$ we have $A *_! B = 0$ if $P_2 \neq P_3$ and $A *_! B = A *_!^{P_2} B$ if $P_2 = P_3$. 
\end{definition} 

\subsection{Rigidity and duality}

Recall that we write 
$$\Dual(-) = \mathcal{H}om(-, \omega),$$
where $\omega = p^!\un$ and $p : \PGQ \to P \backslash \pt/Q$. 

We write $\op{inv} : \PGQ \to Q \backslash G/P$ the map induced by $g \mapsto g^{-1}$ and we write 
$$\Dual^-(-) = \op{inv}^*\Dual(-).$$

\begin{lemma}\label{lem:dual-basic-calculation}
Let $\underline{P}$ be a sequence from $P$ to $Q$ and let $\underline{P}^{\mathrm{op}}$ be its opposite. 
\begin{enumerate}
\item The category $\Hecke_{P,Q}$ is stable under $\Dual$. The functor $\Dual^-$ sends $\Hecke_{P,Q}$ to $\Hecke_{Q,P}$. 
\item We have $\Dual(\BS(\underline{P})) = \BS(\underline{P})$ and $\Dual^-(\BS(\underline{P})) = \BS(\underline{P}^{\mathrm{op}})$.
\item There are natural maps 
\begin{align*}
\Dual(A) *_!^{P_2} \Dual(B) &\to \Dual(A *_*^{P_2} B), \\
\Dual(A *_!^{P_2} B) &\to \Dual(A) *_*^{P_2} \Dual(B), \\
\Dual^-(A *_!^{P_2} B) &\to \Dual^-(B) *_*^{P_2} \Dual^-(A),
\end{align*}
that restrict to isomorphisms on $\Hecke_{P_1,P_2} \times \Hecke_{P_2,P_3}$.
\end{enumerate}
\end{lemma}

\begin{proof}
The same proof as in \cite[Lemma 4.15]{eberhardtUniversalKoszulDuality2025a} holds.
\end{proof}

\begin{lemma}\label{lem:calculation-of-dual}
Let $A \in \DK(P_1 \backslash G/P_2), B \in \DK(P_2 \backslash G/P_3)$ as well as $C \in \DK(P_1 \backslash G/P_3)$. Then, there are natural isomorphisms
$$\Hom(A *_!^{P_2} B, C) = \Hom(A, C *_*^{P_3} \Dual^-(B)) = \Hom(B, \Dual^-(A) *_*^{P_1} C).$$
\end{lemma} 

\begin{proof}
The same proof as in \cite[Lemma 4.16]{eberhardtUniversalKoszulDuality2025a} holds.
\end{proof}

\begin{theorem}\label{thm:Rigitidy-Of-Hecke-Categorioid}
There is a natural equivalence of functors $\Hecke_{P_1,P_2} \times \Hecke_{P_2,P_3} \to \Hecke_{P_1,P_3}$,
$$- *^{P_2}_! - \xrightarrow{\sim} - *_*^{P_2} - .$$
In particular, the total Hecke category $\Hecke_{\tot}$ is rigid with left and right dual canonically isomorphic to $\Dual^-$. 
\end{theorem} 

\begin{proof}
If $G$ is finite dimensional, then the Theorem is clear by  \Cref{lem:calculation-of-dual,lem:comparison-convol-finite-dim}. The total Hecke category $\Hecke_{\tot}$ is generated as a monoidal category (using $*_!$ convolution) by Bott--Samelson objects. By Lemma \ref{lem:dual-basic-calculation} Bott--Samelson objects attached to length one sequences are supported on finite dimensional reductive subgroups of $G$. Therefore these Bott--Samelson objects are dualizable. By Lemma \ref{lem:Basics-of-Bott--Samelson}, any Bott--Samelson is an iterated convolution of length one Bott--Samelson objects. Hence all Bott--Samelson are dualizable. By Lemma \ref{lem:Generateting-families} Bott--Samelson objects generate $\Hecke_{\tot}$ so all objects of $\Hecke_{\tot}$ are dualizable.

Arguing as in \cite[Theorem 4.17]{eberhardtUniversalKoszulDuality2025a}, the category $\Hecke_{\tot}$ is a monoidal category (with $*_!$-convolution) with a dualizing functor $\Dual^-$ in the sense of \cite{boyarchenkoDualityFormalismSpirit2013} and a dualizing object given by $\bigoplus_P \Delta_1^{P,P}$. Such a category comes equipped with two tensor structures which are here identified with $*_!$ and $*_*$ and a canonical map $*_! \to *_*$. Furthermore, by \cite[Lemma 3.4]{boyarchenkoDualityFormalismSpirit2013}, for any object $A \in \Hecke$, the natural transformation $A *_! - \to A *_* -$ (resp. $ -*_! A \to - *_* A$) is an isomorphism if and only if $A$ is left (resp. right) dualizable. As all objects of $\Hecke$ are dualizable, for all objects $A,B$, the map $A *_! B \to A *_* B$ is an isomorphism. 
\end{proof}

\begin{remark}
From now on, we will simply write $*^{P}$ for the convolution when applying it to objects in the total Hecke category.
\end{remark}

\subsection{Purity and formality}
Bott--Samelson objects are pure, in the sense that they arise as the pushforward of the constant object along a proper map with smooth source. 

\begin{remark}\label{rem:weightstructures}
One way to make this notion of purity precise is to construct a Chow weight structure, say $\op{w}_{\op{Chow}}$, in the sense of Bondarko \cite{bondarkoWeightStructuresVs2010} on the category of $K$-motives $\DKbig$. In this way, pure objects correspond to objects in the heart of the weight structure.

The construction of a weight structure on all of $\DKbig$ would rely on deep results on the existence of equivariant resolutions, see \cite{aranhaChowWeightStructure2025} for a similar situation.
However, one might simply reverse engineer the Chow weight structure on the Hecke category. For this, one observes that the Hecke category is glued by categories of the form $\DKT(P\backslash PwQ/Q)\simeq \Perf(\RR_{Q_w})$ for each stratum. The Chow weight structure on the category $\DKT(P\backslash PwQ/Q)$ corresponds to the natural weight structure on $\Perf(\RR_{Q_w})$ whose heart corresponds to projective $\RR_{Q_w}$-modules in cohomological degree $0$. Now one obtains the weight structure $\op{w}_{\op{Chow}}$ on the Hecke category by gluing. One can then show that Bott--Samelson objects generate the heart of this weight structure (this is indeed implied by the stronger point-wise purity result \Cref{thm:pointwisepuritybottsamelson}).

We note that in the context of constructible sheaves, the theory of parity sheaves \cite{juteauParitySheaves2014} is a $\Z/2$-periodic shadow of this idea, for a formal statement in this direction, see \cite[Theorem 4.5]{eberhardtMixedMotivesGeometric2019}.
\end{remark}

In the light of \Cref{rem:weightstructures}, we will simply make the following, most straightforward definition, inspired by the definition of category $\mathcal{K}$ in \cite{soergelRelationIntersectionCohomology2000a}.
\begin{definition}
We define the category of pure objects 
$$\Hecke_{P,Q}^{\pure} = \langle \BS(\underline{P}) \rangle_{\op{add},\op{thick}}$$
to be the full subcategory of $\Hecke_{P,Q}$ generated under direct sums and direct summands by Bott--Samelson sheaves. 
\end{definition}  
\begin{lemma}\label{lem:permanencePropretiesOfPurity}
\begin{enumerate}
\item For $P_1 \subset P_2, Q_1 \subset Q_2$, the functors $^{P_1}_{P_2}\pi_{Q_2}^{Q_1,*}$ and $^{P_1}_{P2}\pi_{Q_2,*}^{Q_1}$ preserve $\Hecke^{\pure}$. 
\item The category $\Hecke^{\pure}_{\tot}$ is closed under convolution. 
\item The category $\Hecke^{\pure}_{P,Q}$ is generated under direct sums and direct summands by objects of the form $^B_P\pi^B_{Q,*}A$ for $A \in \Hecke_{B,B}^{\pure}$.
\end{enumerate}
\end{lemma} 
\begin{proof}
By Lemma \ref{lem:Basics-of-Bott--Samelson}, $^{P_1}_{P_2}\pi_{Q_2}^{Q_1,*}$ and $^{P_1}_{P2}\pi_{Q_2,*}^{Q_1}$ can be realized as convolution against some Bott--Samelson objects, this shows that items $(1)$ and $(2)$ are equivalent. Furthermore, by the same Lemma, the subcategories of Bott--Samelson objects are stable under convolution, this implies item $(2)$. For item $(3)$, it follows from \Cref{cor:projectivebundleformula} that for any $A' \in \Hecke_{P,Q}^{\pure}$, the object $^B_P\pi^B_{Q,*}{\phantom{}^B_P\pi^{B,*}_Q} A'$ is a direct sum of copies of $A'$.
\end{proof}
Essentially the cellular structure of Bott--Samelson resolutions implies that they will fulfill the following pointwise purity condition.
\begin{definition}[Pointwise purity]
An object $A \in \Hecke_{P,Q}$ is $?$-pointwise pure for $? \in \{!,*\}$ if for all $w \in {^PW^Q}$, the object
$i_w^?A$
are direct summands of sums of $\un \in \DK(P \backslash PwQ/Q)$.
\end{definition}
\begin{remark}
	The category of Tate objects $\DKT(P \backslash PwQ/Q)$ is equivalent to the perfect derived category of $A=K_0(P \backslash PwQ/Q)$. Under this equivalence, pointwise pure objects correspond to projective $A$-modules in cohomological degree zero.
\end{remark}
We will make use of the following useful characterization of pointwise purity.
\begin{lemma}[\protect{\cite[Lemma 4.11]{eberhardtUniversalKoszulDuality2025a}}]\label{lem:pointwise-purity-and-delta-flags}
Pointwise $!$-pure objects in $\Hecke$ are exactly the objects that admit a filtration whose associated graded are direct summands of sums of costandard objects. Dually, pointwise $*$-pure objects are exactly the objects that admit a cofiltration whose associated graded are direct summands of sums of standard objects. 
\end{lemma} 
Pointwise purity allows us to prove that mapping spaces are discrete.
\begin{lemma}[\protect{\cite[Lemma 4.12]{eberhardtUniversalKoszulDuality2025a}}]
Let $M,N \in \Hecke_{P,Q}$ be pointwise $*$-pure and pointwise $!$-pure respectively, then 
$\Map(M,N)$
is concentrated in degree $0$. 
\end{lemma}
\begin{proof}
The proof of \cite[Lemma 4.12]{eberhardtUniversalKoszulDuality2025a} holds as the only non-formal input is that the cohomology (in the $6$-formalism $\DK$) of the strata of $\PGQ$ are concentrated in degree $0$. However, we have 
$$\Map_{\DK(P\backslash PwQ/Q)}(\un, \un) = K(P \backslash PwQ/Q)_{\op{r}},$$
and by Lemma \ref{lem:DescriptionOfDKT}, the reduced $K$-theory $K(P \backslash PwQ/Q)_{\op{r}}$ is in degree $0$. 
\end{proof}

\begin{theorem} \label{thm:pointwisepuritybottsamelson}
All objects of $\Hecke_{P,Q}^{\pure}$ are pointwise $!$ and $*$-pure. 
\end{theorem}

\begin{proof}
For $P = Q = B$, this is done \cite[Lemma 4.13]{eberhardtUniversalKoszulDuality2025a} for $\Hecke_{B,B}^{\circ}$. The same property holds on all connected components of $G$ by Lemma \ref{lem:Action-Of-Min-Elts}. In view of $(3)$ of Lemma \ref{lem:permanencePropretiesOfPurity}, it is enough to prove that the functor $^B_P\pi^B_{Q,*}$ preserves $!$ and $*$-pointwise pure objects. Let us shorten $^B_P\pi^B_{Q} = \pi$. By Lemma \ref{lem:ChangeOfParabolicPushforward}, $\pi_*$ sends (co)standards to direct summands of sums of (co)standards and therefore sends $\Delta$ (resp. $\nabla$)-filtered objects to $\Delta$ (resp. $\nabla$)-filtered objects. The statement then follows from the description in Lemma \ref{lem:pointwise-purity-and-delta-flags}.
\end{proof}

We equip the category $\Hecke_{P,Q}$ with the weight structure constructed as in \cite[Proposition A.4]{eberhardtUniversalKoszulDuality2025a}, its heart is $\Hecke_{P,Q}^{\pure}$ and as in \cite[Corollary 4.14]{eberhardtUniversalKoszulDuality2025a} we have: 

\begin{theorem}\label{corol:weight-complex-functor-is-an-equiv}
The weight complex functor yields an equivalence of categories
$$\Hecke_{P,Q} \xrightarrow{\sim} \Ch^b(\Hecke_{P,Q}^{\pure}).$$
\end{theorem}

\subsection{$\Kyp$-functors}
Our next goal is to compare the pure objects in the singular $K$-motivic Hecke category to singular $K$-theory Soergel bimodules. For this, we use an analog, denoted by $\Kyp$, to the `hypercohomology' functor $\mathbb{H}$ as well as Soergel's Erweiterungssatz, see \cite{soergelKategorieMathcalPerverse1990}. Our first main goal is to investigate the monoidal structure of this functor.
\begin{definition}
The $\Kyp$-functor of type $(P,Q, \gamma)$ is the functor 
$$\Hecke_{P,Q}^{\gamma} \to \Dd(\RR_P \otimes \RR_Q)$$ 
given by 
$$\Kyp_{P,Q}^{\gamma}(A) = \Map(\un, \alpha_{P,Q,*}^{\gamma}A)$$
where $\alpha_{P,Q}^{\gamma} : P \backslash G^{\gamma} /Q \to P \backslash\pt/Q$ is the natural structure map. 
\end{definition} 

Similarly, we can globalize the various $\Kyp$-functors as
$$\bigoplus_{P,Q,\gamma} \Hecke_{P,Q}^{\gamma} \xrightarrow{\oplus\Kyp_{P,Q}^{\gamma}} \bigoplus_{P,Q,\gamma} \Dd(\RR_P \otimes \RR_Q).$$
We equip the category $\bigoplus_{P,Q,\gamma} \Dd(\RR_P \otimes \RR_Q)$ with the `$\Gamma$-graded convolution of bimodules monoidal structure'. Let us write $\Dd(\RR_P \otimes \RR_Q)^{\gamma}$ to be the summand indexed by $\gamma \in \Gamma$. Then the convolution of two objects $A_1 \in \Dd(\RR_{P_1} \otimes \RR_{P_2})^{\gamma_1}$ and $A_2 \in \Dd(\RR_{P_3} \otimes \RR_{P_4})^{\gamma_2}$ is 
\begin{enumerate}
\item $A_1 * A_2 = 0$ if $P_2 \neq P_3$, 
\item $A_1 * A_2 = A_1 \otimes_{\RR_{P_2}} A_2 \in \Dd(\RR_{P_1} \otimes \RR_{P_4})^{\gamma_1\gamma_2}$ otherwise. 
\end{enumerate}

\begin{lemma}\label{lem:FactoThroughDKT}
For all $A \in \Hecke_{P,Q}^{\gamma}$, the object $\alpha_{P,Q,*}^{\gamma}A \in \DK(\pt/P \times \pt/Q)$ lies in $\DKT(\pt/P \times \pt/Q)$. 
\end{lemma}

\begin{proof}
It is enough to prove that $\alpha_{P,Q,*}^{\gamma}\Delta_w^{P,Q}$ lies in $\DKT$. By Corollary \ref{corol:Generation-By-BB} and Lemma \ref{lem:FunctorialityOfDKT}, we can assume that $P = Q = B$. In this case, the Bruhat strata are affine spaces and the Lemma follows.
\end{proof}

We start by constructing a lax monoidal structure on the functor $\Kyp$. Let us fix $P_1, P_2$ and $P_3$ three parabolics and $\gamma_1, \gamma_2 \in \Gamma$. We will shorten $\alpha_{P_1,P_2}^{\gamma_1} = \alpha_1, \alpha_{P_2,P_3}^{\gamma_2} = \alpha_2$ and $\alpha_{P_1,P_3}^{\gamma_1\gamma_2} = \alpha_3$. Let us write the following maps 
\[\begin{tikzcd}
	{P_1 \backslash G^{\gamma_1\gamma_2} /P_3} & {P_1 \backslash \pt/P_3} \\
	{P_1 \backslash G^{\gamma_1} \times^{P_2} G^{\gamma_2} /P_3} & {P_1 \backslash \pt \times^{P_2} \pt /P_3} \\
	{P_1 \backslash G^{\gamma_1} /P_2 \times P_2 \backslash  G^{\gamma_2} /P_3} & {P_1 \backslash \pt /P_2 \times P_2 \backslash \pt /P_3}
	\arrow["{\alpha_3}", from=1-1, to=1-2]
	\arrow["m", from=2-1, to=1-1]
	\arrow["{\tilde{\alpha}}"', from=2-1, to=2-2]
	\arrow["{\pr_1\times \pr_2}"', from=2-1, to=3-1]
	\arrow["{\tilde{m}}"', from=2-2, to=1-2]
	\arrow["{\tilde{\pr}}", from=2-2, to=3-2]
	\arrow["{\alpha_1 \times \alpha_2}"', from=3-1, to=3-2]
\end{tikzcd}\]
Note that the bottom square is Cartesian. To construct the lax monoidal structure on $\Kyp$, we equip $\oplus_{P,P', \gamma} \alpha_{P,P'}^{\gamma,*}\un$ with a natural coalgebra structure. Note that, unless $G$ is finite dimensional, this object does not belong to $\Hecke$, but it still lies in $\oplus_{P,P'} \DK(P \backslash G/P')$. Specifically, we construct a map
\begin{equation}\label{eq:coalg-structure}
\alpha_3^*\un \to \alpha_1^*\un *_*^{P_2} \alpha_{2}^*\un.
\end{equation}
Observe that the RHS is canonically 
$$\alpha_1^*\un *_*^{P_2} \alpha_{2}^*\un = m_*(\pr_1 \times \pr_2)^!(\alpha_1 \times \alpha_2)^*\un$$
by adjunction, it is enough to construct a natural map 
$$(\pr_1 \times \pr_2)_!m^*\alpha_3^*\un = (\pr_1 \times \pr_2)_!\un \to (\alpha_1 \times \alpha_2)^*\un = \un.$$
This is the map dual to 
$$\un \to (\pr_1 \times \pr_2)_*\un.$$
It is an elementary verification that this defines a canonical coalgebra structure on $\oplus_{P,P', \gamma} \alpha_{P,P'}^{\gamma,*}\un$.

Let us also shorten $\Kyp_{P_1,P_2}^{\gamma_1} = \Kyp_1, \Kyp_{P_2, P_3}^{\gamma_2} = \Kyp_2$ and $\Kyp_{P_1,P_3}^{\gamma_1\gamma_2} = \Kyp_3$. The lax-monoidal structure is then given by the following map 
\begin{align*}
\Kyp_1(-) \otimes \Kyp_2(-) &= \Map(\un, \alpha_{1,*}-) \otimes \Map(\un, \alpha_{2,*}-) \\
&\simeq \Map(\alpha_1^*\un, -) \otimes \Map(\alpha_2^*\un, -) \\
&\to \Map(\alpha_1^*\un *_{*}^{P_2} \alpha_2^*\un, - *^{P_2} -) \\
&\to \Map(\alpha_3^*\un, - *^{P_2} -) \\
&= \Kyp_3(- *^{P_2} -),
\end{align*}
where the first map is given by adjunction, the second is the functoriality of convolution, the third is precomposition by \eqref{eq:coalg-structure}. 

\begin{lemma}\label{lem:calculation-for-bott-samelson}
We keep the notation as above and the chosen $P_1,P_2,P_3$ and $\gamma_1,\gamma_2$. Let us choose $\underline{Q}$ and $\underline{Q}'$ Bott--Samelson sequences as before from $P_1 \to P_2$ and $P_2 \to P_3$ sitting above $\gamma_1$ and $\gamma_2$ respectively. There is a canonical isomorphism
$$\Kyp_1(\BS(\underline{Q})) \otimes_{\RR_{P_2}} \Kyp_2(\BS(\underline{Q}')) \simeq \Kyp_3(\BS(\underline{Q} \sqcup \underline{Q}'))$$
making the following diagram commutative 
\[\begin{tikzcd}
	{\Kyp_1(\BS(\underline{Q})) \otimes \Kyp_2(\BS(\underline{Q}'))} & {\Kyp_3(\BS(\underline{Q}) *^{P_2} \BS(\underline{Q}'))} \\
	{\Kyp_1(\BS(\underline{Q})) \otimes_{\RR_{P_2}} \Kyp_2(\BS(\underline{Q}'))} & { \Kyp_3(\BS(\underline{Q} \sqcup \underline{Q}'))}
	\arrow["\lax", from=1-1, to=1-2]
	\arrow[from=1-1, to=2-1]
	\arrow[equals, from=1-2, to=2-2]
	\arrow["\sim"', from=2-1, to=2-2]
\end{tikzcd}\]
where the top vertical map is the lax monoidal structure on $\Kyp$ previously constructed and the left vertical map is the natural map $\otimes \to \otimes_{\RR_{P_2}}$. 
\end{lemma} 

\begin{proof}
Let us start by considering the following diagram 
\[\begin{tikzcd}[column sep = 20 pt]
	{P_1 \backslash \underline{Q \sqcup Q}'/P_3} & {P_1 \backslash G^{\gamma_1\gamma_2} /P_3} & {P_1 \backslash \pt/P_3} \\
	{P_1 \backslash \underline{Q} \times^{P_2} \underline{Q}'/P_3} & {P_1 \backslash G^{\gamma_1} \times^{P_2} G^{\gamma_2} /P_3} & {P_1 \backslash \pt \times^{P_2} \pt /P_3} \\
	{P_1 \backslash \underline{Q}/P_2 \times P_2 \backslash \underline{Q}'/P_3} & {P_1 \backslash G^{\gamma_1} /P_2 \times P_2 \backslash  G^{\gamma_2} /P_3} & {P_1 \backslash \pt /P_2 \times P_2 \backslash \pt /P_3}
	\arrow["{\pi_{Q \sqcup Q'}}", from=1-1, to=1-2]
	\arrow["{\alpha_3}", from=1-2, to=1-3]
	\arrow["m", from=2-1, to=1-1]
	\arrow["{\pi_{Q,Q'}}", from=2-1, to=2-2]
	\arrow["{\pr_{Q,Q'}}"', from=2-1, to=3-1]
	\arrow["m", from=2-2, to=1-2]
	\arrow["{\tilde{\alpha}}"', from=2-2, to=2-3]
	\arrow["{\bar\pr}"', from=2-2, to=3-2]
	\arrow["{\tilde{m}}"', from=2-3, to=1-3]
	\arrow["{\tilde{\pr}}", from=2-3, to=3-3]
	\arrow["{\bar\pi}", from=3-1, to=3-2]
	\arrow["{\bar\alpha}", from=3-2, to=3-3]
\end{tikzcd}\]
where the maps are as defined above. We abbreviate $\bar\pi =\pi_Q\times \pi_{Q'}$, $\bar\alpha=\alpha_1\times \alpha_2$ and $\bar\pr=\pr_1\times \pr_2$ so that the following diagrams fit within a page. Note that the two bottom squares in the diagram are Cartesian. 

We first observe that 
$$\Kyp_1(\BS(\underline{Q})) \otimes \Kyp_2(\BS(\underline{Q}')) = \Map(\un, \bar\alpha_*\bar\pi_*\un).$$
Although the map $\tilde{m}$ is not representable, recall from Section \ref{sec:non-representable-functors} that we defined a pushforward $\tilde{m}_*$ and a pullback $\tilde{m}^*$ on Tate objects. By Lemma \ref{lem:forget-along-diag}, there is a canonical isomorphism 
$$\Map(\un, \tilde{m}_*\tilde{\pr}^*\bar\alpha_*(\pi_{Q} \times \pi_{Q'})_*\un) = \Kyp_1(\BS(\underline{Q})) \otimes_{\RR_{P_2}} \Kyp_2(\BS(\underline{Q}')).$$
We then have
\begin{align*}
 \Kyp_1(\BS(\underline{Q})) \otimes_{\RR_{P_2}} \Kyp_2(\BS(\underline{Q}')) &= \Map(\un, \tilde{m}_*\tilde{\pr}^*\bar\alpha_*\bar\pi_*\un) \\
&\simeq \Map(\un, \tilde{m}_*\tilde{\alpha}_*\pi_{Q,Q',*}\un) \\
&\simeq \Map(\un, \alpha_{3,*}m_*\pi_{Q,Q',*}\un) \\
&= \Kyp_3(\BS(\underline{Q} \sqcup \underline{Q})),
\end{align*}
where the maps are given by 
\begin{enumerate}
\item the first one is the previous observation, 
\item the second is by base change in the above diagram 
\item the third one follows from Lemma \ref{lem:Compatibility-Pullback} and Remark \ref{remark:Compatibility-pushforward} as both maps $\tilde{\alpha}$ and $\alpha_3 \circ m$ are representable.
\end{enumerate}

To prove the compatibility with the lax monoidal structure, we construct the following commutative diagram.

\[\begin{tikzcd}[column sep=60pt]
	{\Kyp_1(\BS(\underline{Q})) \otimes \Kyp_2(\BS(\underline{Q}'))} & {\Kyp_1(\BS(\underline{Q})) \otimes_{\RR_{P_2}} \Kyp_2(\BS(\underline{Q}'))} \\
	{\Map(\un, \bar{\alpha}_*\bar{\pi}_*\un)} & {\Map(\tilde{\pr}_!\tilde{m}^*\un, \bar{\alpha}_*\bar{\pi}_*\un)} \\
	{\Map(\bar{\alpha}^*\un, \bar{\pi}_*\un)} & {\Map(\bar{\alpha}^*\tilde{\pr}_!\tilde{m}^*\un, \bar{\pi}_*\un)} \\
	& {\Map(\bar{\pr}_!\tilde{\alpha}^*\tilde{m}^*\un, \bar{\pi}_*\un)} \\
	{\Map(\bar{\pr}^!\bar{\alpha}^*\un, \bar{\pr}^!\bar{\pi}_*\un)} & {\Map(\tilde{\alpha}^*\tilde{m}^*\un, \bar{\pr}^!\bar{\pi}_*\un)} \\
	{\Map(\bar{\pr}^!\bar{\alpha}^*\un, \pi_{Q,Q',*}\pr_{Q,Q'}^!\un)} & {\Map(\tilde{\alpha}^*\tilde{m}^*\un, \pi_{Q,Q',*}\pr_{Q,Q'}^!\un)} \\
	{\Map(\bar{\pr}^!\bar{\alpha}^*\un, \pi_{Q,Q',*}\un)} & {\Map(\tilde{\alpha}^*\tilde{m}^*\un, \pi_{Q,Q',*}\un)} \\
	{\Map(\bar{\pr}^!\bar{\alpha}^*\un, \pi_{Q,Q',*}\un)} & {\Map(m^*\alpha_3^*\un, \pi_{Q,Q',*}\un)} \\
	{\Map(m_*\bar{\pr}^!\bar{\alpha}^*\un, m_*\pi_{Q,Q',*}\un)} & {\Map(\alpha_{3}^*\un, m_*\pi_{Q,Q',*}\un)} \\
	& {\Kyp_3(\BS(\underline{Q} \sqcup \underline{Q}'))}
	\arrow["{{\otimes \to \otimes_{\RR_{P_2}}}}", from=1-1, to=1-2]
	\arrow[from=1-1, to=2-1]
	\arrow[from=1-2, to=2-2]
	\arrow["{{\tilde{\pr}_!\tilde{m}^*\un \to \un}}", from=2-1, to=2-2]
	\arrow["\adj"', from=2-1, to=3-1]
	\arrow["\adj", from=2-2, to=3-2]
	\arrow["{{\tilde{\pr}_!\tilde{m}^*\un \to \un}}", from=3-1, to=3-2]
	\arrow["{{\bar{\pr}_! \tilde{\alpha}^*\tilde{m}^*\un \to \bar{\alpha}^*\un}}"', from=3-1, to=4-2]
	\arrow["{{\bar{\pr}^!}}"', from=3-1, to=5-1]
	\arrow["\BC", from=3-2, to=4-2]
	\arrow["\adj", from=4-2, to=5-2]
	\arrow["{{\tilde{\alpha}^*\tilde{m}^*\un \to \bar{\pr}^!\bar{\alpha}^*\un}}", from=5-1, to=5-2]
	\arrow["\BC"', from=5-1, to=6-1]
	\arrow["\BC", from=5-2, to=6-2]
	\arrow["{{\tilde{\alpha}^*\tilde{m}^*\un \to \bar{\pr}^!\bar{\alpha}^*\un}}", from=6-1, to=6-2]
	\arrow["{{\pr_{Q,Q'}^!\un = \un}}"', equals, from=6-1, to=7-1]
	\arrow["{{\pr_{Q,Q'}^!\un = \un}}", equals, from=6-2, to=7-2]
	\arrow["{{\tilde{\alpha}^*\tilde{m}^*\un \to \bar{\pr}^!\bar{\alpha}^*\un}}", from=7-1, to=7-2]
	\arrow[equals, from=7-1, to=8-1]
	\arrow["{m^*\alpha_3^*\un = \tilde{\alpha}^*\tilde{m}^*\un}", from=7-2, to=8-2]
	\arrow["{{m^*\alpha_3^*\un \to \bar{\pr}^!\bar{\alpha}^*\un}}"', from=8-1, to=8-2]
	\arrow["{{m_*}}"', from=8-1, to=9-1]
	\arrow["\adj", from=8-2, to=9-2]
	\arrow["{{\alpha_{3}^*\un \to m_*\bar{\pr}^!\bar{\alpha}^*\un}}"', from=9-1, to=9-2]
	\arrow[equals, from=9-2, to=10-2]
\end{tikzcd}\]
Let us first describe the vertical maps. 
The maps labeled $\BC$ are induced by proper base change in one of the diagrams above. Next, the maps labeled $\adj$ are obtained by using some adjunction. Lastly, the maps $\pr_{Q,Q'}^!\un = \un$ and $m^*\alpha_3^*\un = \tilde{\alpha}^*\tilde{m}^*\un$ are simply replacing the LHS with the RHS. 

The horizontal maps are all induced by some precomposition in the first entry in $\Map$ by the map that labels the arrow. We will describe these maps below. Before proving the commutativity of this diagram, observe that the composition along the right vertical maps yields (with extra steps) the isomorphism of the first half of the Lemma and the composition along the left vertical and last bottom will be the lax-monoidality. 

We now describe the remaining maps. The map
\begin{equation}\label{eq:map1}
\tilde{\pr}_!\tilde{m}^*\un \to \un
\end{equation}
is the map adjoint to $\un = \tilde{m}^*\un \to \tilde{\pr}^!\un = \un$. Next, the map 
\begin{equation}\label{eq:map2}
\bar\pr_! \tilde{\alpha}^*\tilde{m}^*\un \to \bar\alpha^*\un 
\end{equation}
is the map obtained as
\begin{align*}
\bar\pr_! \tilde{\alpha}^*\tilde{m}^*\un \xrightarrow{\BC} \bar\alpha^*\tilde{\pr}_!\tilde{m}^*\un \xrightarrow{\bar\alpha^*(\ref{eq:map1})} \bar\alpha^*\un.
\end{align*}
Moreover, we obtain the map 
\begin{equation}\label{eq:map3}
\tilde{\alpha}^*\tilde{m}^*\un \to \bar\pr^!\bar\alpha^*\un
\end{equation}
as adjoint to the map \eqref{eq:map2}. The map 
\begin{equation}\label{eq:map4}
m^*\alpha_3^*\un \to \bar{\pr}^!\bar{\alpha}^*\un
\end{equation} 
is obtained by using that $\tilde{\alpha}^*\tilde{m}^*\un = m^*\alpha_{3}^*\un = m^*\un$ and using  \eqref{eq:map3}. The map 
\begin{equation}\label{eq:map5}
{\alpha_{3}^*\un \to m_*\bar{\pr}^!\bar{\alpha}^*\un}
\end{equation}
is the map adjoint to \eqref{eq:map4}.
\end{proof}

\begin{theorem}\label{thm:Strict-Monoidality}
The functor $\Kyp$ is strict monoidal. 
\end{theorem}

\begin{proof}
This is an immediate consequence of Lemma \ref{lem:calculation-for-bott-samelson}. 
\end{proof}

\begin{remark}\label{remark:Image-Of-K-is-SBim}
It also follows from Lemma \ref{lem:calculation-for-bott-samelson} that $\Kyp(\BS(\underline{Q}))$ is a Bott--Samelson bimodule.
\end{remark}

\begin{remark}\label{remark:compatibility-with-pullpush}
It follows from Remark \ref{remark:Example-Bott--Samelson} and Lemma \ref{lem:calculation-for-bott-samelson}, that if we let $\pi = {\phantom{}^P_{P'}}\pi^Q_{Q'} : P \backslash G/Q \to P' \backslash G/Q'$, then there are canonical isomorphisms (for all $\gamma \in \Gamma$): 
\begin{enumerate}
\item $\Kyp_{P,Q}^{\gamma}(\pi^*-) = \RR_{P} \otimes_{\RR_{P'}} \Kyp_{P',Q'}^{\gamma}(-) \otimes_{\RR_{Q'}} \RR_Q$, 
\item $\Kyp_{P',Q'}^{\gamma}(\pr_*-) = \Kyp_{P,Q}^{\gamma}(-)$ as $\RR_{P'} \otimes \RR_{Q'}$-modules.
\end{enumerate}
\end{remark} 

\begin{construction}
Since the weight complex functor is monoidal by \cite{aokiWeightComplexFunctor2020}, it induces a natural transformation of functors $\Parab(G) \to \Morita(\Pr^{\omega})$
$$\Hecke \Rightarrow \Ch^b(\Hecke^{\pure}).$$
Since the functor $\Kyp$ is strictly monoidal by Theorem \ref{thm:Strict-Monoidality} and sends pure objects to Soergel bimodules by Remark \ref{remark:Image-Of-K-is-SBim}, it induces a natural transformation of functors $\Parab(G) \to \Morita(\Pr^{\omega})$, 
$$\Ch^b(\Hecke^{\pure}) \Rightarrow \Ch^b(\SBim).$$
Composing these two functors, we get a natural transformation, still denoted by $\Kyp$ 
\begin{equation}\label{eq:defi-K-functors}
\Hecke \Rightarrow \Ch^b(\SBim).
\end{equation}
\end{construction}

\subsection{Equivalence of categories} We are now ready to show Soergel's Erweiterungssatz in our setting.
\begin{theorem}\label{thm:K-is-an-equivalence}
The natural transformation \eqref{eq:defi-K-functors} 
$$\Kyp : \Hecke \Rightarrow \Ch^b(\SBim)$$
is an isomorphism. 
\end{theorem}

\begin{remark}
The statement of Theorem \ref{thm:K-is-an-equivalence} amounts to proving that for all $P,Q$ and $\gamma \in \Gamma$, the functor
$$\Kyp_{P,Q}^{\gamma} : \Hecke_{P,Q}^{\gamma} \rightarrow \Ch^b(\SBim_{P,Q}^{\gamma})$$
is an equivalence. 
\end{remark}

\begin{theorem}\label{thm:Fully-Faithfulness}
For all $P,Q, \gamma$, the functor $\Kyp_{P,Q}^{\gamma} : \Hecke_{P,Q}^{\gamma, \pure} \to \SBim_{P,Q}^{\gamma}$ is fully faithful. 
\end{theorem}

\begin{proof}[Proof of Theorem \ref{thm:Fully-Faithfulness}]
We will use an argument very similar to \cite[Proposition III.6.11]{soergelEquivariantMotivesGeometric2018} (`geklaut und umgespritzt'). 
	Let $M,N\in \Hecke_{P,Q}^{\gamma, \pure}.$ Using $M=\Delta^{P,P}_e*^P M$, $R_P=\Kyp_{P,P}^{e}(\Delta^{P,P}_e)$ as well as Theorem \ref{thm:Strict-Monoidality}, we obtain the following commutative diagram
	\[\begin{tikzcd}
	{\Hom_{\Hecke_{P,Q}}(M, M')} & {\Hom_{R_P\otimes R_Q}(\Kyp_{P,Q}^{\gamma}(M), \Kyp_{P,Q}^{\gamma}(M'))} \\
	{\Hom_{\Hecke_{P,P}}(\Delta^{P,P}_e,N) } & {\Hom_{R_P\otimes R_P}(R_P,\Kyp_{P,P}^{\circ}(N)))}
	\arrow["\Kyp", from=1-1, to=1-2]
	\arrow["\wr"', from=1-1, to=2-1]
	\arrow["\wr", from=1-2, to=2-2]
	\arrow["\Kyp", from=2-1, to=2-2]
\end{tikzcd}\]
	where $N=M*^P\Dual^-(M')\in\Hecke_{P,P}^\circ.$ Hence, we reduced the statement to showing that the bottom vertical arrow is an equivalence. Since this only involves the identity component, we can assume that $G$ is connected for now. Moreover, we have shown how to reduce from general $P,Q$ to $P=Q$.
	In the case of $P=B$, the statement is \cite[Theorem 4.22]{eberhardtUniversalKoszulDuality2025a}. We will explain how to reduce to this case using the projective bundle formula. For this, we abbreviate $\pi: B\bs G/P\to P\bs G/P$ which is a $P/B$-bundle. By \Cref{cor:projectivebundleformula}, $\pi_*\pi^*\Delta_e^{P,P}=(\Delta_e^{P,P})^{\oplus n}$ for $n=|W/W_P|.$ Then, by Remark \ref{remark:compatibility-with-pullpush}, we obtain the commutative diagram
	
\[\begin{tikzcd}
	{\Hom_{\Hecke_{B,P}}(\pi^*\Delta^{P,P}_e,\pi^*N)} & {\Hom_{R_B\otimes R_P}(\Kyp_{B,P}(\pi^*\Delta^{P,P}_e),\Kyp_{B,P}(\pi^*N))} \\
	{\Hom_{\Hecke_{P,P}}(\pi_*\pi^*\Delta^{P,P}_e,N)} & {\Hom_{R_P\otimes R_P}(\Kyp_{P,P}(\pi_*\pi^*\Delta^{P,P}_e),N)} \\
	{\Hom_{\Hecke_{P,P}}(\Delta^{P,P}_e,N)^{\oplus n}} & {\Hom_{R_P\otimes R_P}(\Kyp_{P,P}(\Delta^{P,P}_e),N)^{\oplus n}}
	\arrow[from=1-1, to=1-2]
	\arrow["\wr", from=1-1, to=2-1]
	\arrow["\wr", from=1-2, to=2-2]
	\arrow[from=2-1, to=2-2]
	\arrow["\wr", from=2-1, to=3-1]
	\arrow["\wr", from=2-2, to=3-2]
	\arrow[from=3-1, to=3-2]
\end{tikzcd}\]
where the bottom map applies $\Kyp_{P,P}$ to each component. Since each component of the map is an isomorphism if the direct sum is, we reduced to showing that the top horizontal map is an isomorphism. So it suffices to show that $\Kyp_{B,P}$ is fully faithful. By the argument in step one, it suffices to show that $\Kyp_{B,B}$ is fully faithful and we are done.
\end{proof}
\begin{proof}[Proof of Theorem \ref{thm:K-is-an-equivalence}]
We first argue that for all $P,Q,\gamma$, the functor $\Kyp_{P,Q}^{\gamma} : \Hecke_{P,Q}^{\gamma} \rightarrow \SBim_{P,Q}^{\gamma}$ is an equivalence. By Theorem \ref{thm:Fully-Faithfulness}, it is fully faithful. It follows that $\Kyp_{P,Q}^{\gamma}$ yields an equivalence on subcategories of $\Hecke_{P,Q}^{\gamma}$ and $\SBim_{P,Q}^{\gamma}$ of Bott--Samelson bimodules. The claim follows from the fact that both $\Hecke_{P,Q}^{\gamma}$ and $\SBim_{P,Q}^{\gamma}$ are the Karoubian completions of the subcategories of Bott--Samelson bimodules. 

It now follows that the functor $\Ch^b(\Hecke_{P,Q}^{\gamma, \pure}) \rightarrow \Ch^b(\SBim_{P,Q}^{\gamma})$ is an equivalence. Finally, since the weight complex functor $\Hecke_{P,Q}^{b} \rightarrow \Ch^b(\Hecke_{P,Q}^{\gamma, \pure})$ is an equivalence by Corollary \ref{corol:weight-complex-functor-is-an-equiv}, the Theorem follows. 
\end{proof}